\documentclass[a4paper,oneside,11pt]{amsart}
\usepackage{enumerate,stackrel}

\usepackage[english]{babel}
\usepackage{fancyhdr}
\usepackage{amsmath}
\usepackage{mathrsfs}
\usepackage{color}

\usepackage{comment}
\let\oldmarginpar\marginpar
\renewcommand\marginpar[1]{\-\oldmarginpar[\raggedleft\small\sf
#1] {\raggedright\small\sf #1}}

\usepackage{amssymb}

\newtheorem{definition}{Definition}[section]

\newtheorem{theorem}[definition]{Theorem}

\newtheorem{corollary}[definition]{Corollary}
\newtheorem{lemma}[definition]{Lemma}
\newtheorem{remark}[definition]{Remark}

\numberwithin{equation}{section}

\def\bkR{{\rm I\kern-.17em R}}

\def\bkN{{\rm I\kern-.17em N}}

\def\bkC{{\rm \kern.24em \vrule width.05em height1.4ex depth-.05ex
     \kern-.26em C}}

\def\bkK{{\rm I\kern-.22em K}}
\def\bkP{{\rm I\kern-.22em P}}

\newcommand{\mcal}{\mathcal}

\newcommand{\bpr}{\begin{proof}}
\newcommand{\epr}{\hfill $\square$ \end{proof}}
\newcommand{\bprof}{\begin{proofof}}
\newcommand{\eprof}{\hfill $\square$ \end{proofof}}
\newcommand{\bproff}{\begin{proofoff}}
\newcommand{\eproff}{\hfill $\square$ \end{proofoff}}
\begin{document}


\title{Relationship between limiting K-spaces and J-spaces in the real interpolation} 


\author{Bohum\'{\i}r Opic}
\author{Manvi Grover}


\address{Bohum\'{\i}r Opic, 
Department of Mathematical Analysis, 
Faculty of Mathematics and Physics, Charles University,
 Sokolovsk\'a 83, 186 75 Prague 8, Czech Republic}
\email{opic@karlin.mff.cuni.cz
\vskip0,08cm 
}

\address{Manvi Grover, 
Department of Mathematical Analysis, 
Faculty of Mathematics and Physics, Charles University,
 Sokolovsk\'a 83, 186 75 Prague 8, Czech Republic}
\email{grover.manvi94@gmail.com}



\begin{abstract}
In the paper Description of the $K$-spaces by means of $J$-spaces and the reverse problem, Math. Nachr. 296 (2023), no. 9, 4002--4031,  we have establish conditions under which the limiting $K$-space $(X_0,X_1)_{0,q,b;K}$, involving  
a slowly varying function $b$, can be described by means of the $J$-space $(X_0,X_1)_{0,q,a;J}$, with a convenient slowly varying function $a$, and we have also solved the reverse problem. It has been shown that if these conditions are not satisfied that the given problem may not have a solution.

In this paper we assume that these conditions are not satisfied. Nevertheless, our aim is to express the limiting $K$-space $(X_0,X_1)_{0,q,b;K}$ 
as some limiting $J$-space $(Y_0,Y_1)_{0,q,A;J}$, and, similarly, to express the limiting $J$-space $(X_0,X_1)_{0,q,a;J}$ as
 a convenient limiting $K$-space $(Z_0,Z_1)_{0,q,B;K}$. To be more precise, we show that 
$$(X_0,X_1)_{0,q,b;K}= (X_0,X_0+X_1)_{0,q,A;J}=X_0+(X_0,X_1)_{0,q,A;J}$$
and 
$$(X_0,X_1)_{0,q,a;J}= (X_0,X_0\cap X_1)_{0,q,B;K}=X_0\cap (X_0,X_1)_{0,q,B;K},$$
where $A$ and $B$ are convenient weights. Moreover, we establish equivalent norms in the above mentioned spaces. 
The obtained results are applied to get density theorems for spaces in question.



\end{abstract}

\thanks{The research has been supported by the grant no. 23-04720S of the Czech Science Foundation.} 

\maketitle


\section{Introduction and main results}\label{section1}

The real interpolation method $(X_0,X_1)_{\theta,q}$ plays a very important role in numerous 
parts of mathematics 
(cf., e.g., monographs  \cite{BB}, \cite{BL}, \cite{Tri78}, \cite{KPS78}, \cite{BS:IO}, \cite{BK}, \cite{L18}). The definitions of $(X_0,X_1)_{\theta,q}$ by means of 
the Peetre $K$- and $J$-functionals are the most familiar. If  
 $0<\theta< 1$ and $1\le q\le\infty$, then, by the equivalence theorem (see, e.g.,  \cite [Chapter 5, Theorem 2.8]{BS:IO}), these two constructions give the same spaces. 

However, some problems in mathematical analysis have motivated the investigation of the real interpolation with the limiting values $\theta=0$ or $\theta=1$.  
Nowadays, there is a~lot of works, where the limiting real interpolation is studied or applied (see, e.g., \cite{GM86}, \cite{Mil94}, \cite{Dok91}, \cite{EO00:RILFR}, \cite{EOP02:RILF}, \cite{GOT2002Lrrinwsvf}, \cite{CFCKU}, \cite{CFCM10}, \cite{AEEK11}, \cite{CK}, 
\cite{FMS12}, \cite{EO14}, \cite{FMS14}, 
\cite{COBOS201543}, \cite{CS}, \cite{CFCM15}, \cite{FMS15}, \cite{CD3}, \cite{Dok18}, \cite{FFGKR18}, \cite{CFC19}, \cite{AFH20}, 
\cite{BCFC20}, \cite{NO20}, \cite{ODST2019MAMS}, \cite{OG}, \cite{CFCG24}, and the references given there).

In order that the spaces $(X_0,X_1)_{\theta,q}$ are meaningful also with limiting values $\theta=0$ or $\theta=1$ in a general case, one has to extend the given interpolation functor by a convenient weight $v$. This weight $v$ usually belongs to the class $SV(0,\infty)$ of slowly varying functions on the interval $(0,\infty)$. If $\theta \in (0,1)$, then the corresponding space $(X_0,X_1)_{\theta,q,v}$ is a particular case of an interpolation space 
with a function parameter and, by \cite[Theorem 2.2]{Gu78}, the equivalence theorem continues to hold. However, if $\theta=0$ or $\theta=1$, then the corresponding space depends on the fact whether the $K$- or $J$- functional is used to define it; accordingly, the resulting space is denoted by  
$(X_0,X_1)_{\theta,q,v;K}$ or by $(X_0,X_1)_{\theta,q,v;J}$. 

Now a natural question arises: Given  $\theta \in \{0,1\}, q \in [1,\infty]$ and $v \in SV(0,\infty)$, can we describe the space $(X_0,X_1)_{\theta,q,v;K}$ as a $(X_0,X_1)_{\theta,q,w;J}$ space with a convenient $w \in SV(0,\infty)$?



Note that it is sufficient to study the given question only for $\theta=0$ since the answer for $\theta=1$ follows from the solution with $\theta=0$ and from the equality  
$(X_0,X_1)_{0,q,v;K}=(X_1,X_0)_{1,q,u;K}$, or $(X_0,X_1)_{0,q,v;J}=(X_1,X_0)_{1,q,u;J}$,  
 where $u(t)=v(1/t)$ for all $t>0$ (which is a consequence of the fact that $K(f,t;X_0,X_1)=tK(f,t^{-1};X_1,X_0)$ if  
$f \in X_0+X_1$ and $t>0$, or 
$J(f,t;X_0,X_1)=tJ(f,t^{-1};X_1,X_0)$ if $f \in X_0\cap X_1$ and $t>0$, and a change of variables).
 
The partial answer to the given question has been given in \cite{OG}, where we have establish conditions under which the limiting $K$-space $(X_0,X_1)_{0,q,b;K}$, involving  a slowly varying function $b$, can be described by means of the $J$-space $(X_0,X_1)_{0,q,a;J}$, with a convenient slowly varying function $a$. It has been also shown that if these conditions are not satisfied, then the given problem may not have a solution. Moreover, in \cite{OG} also the reverse problem has been solved, i.e., we have establish conditions under which 
the limiting $J$-space $(X_0,X_1)_{0,q,a;J}$, involving a slowly varying function $a$, can be described by means of the 
$K$-space $(X_0,X_1)_{0,q,b;K}$, with a convenient slowly varying function $b$. Then we have applied our results to obtain density theorems for the corresponding limiting interpolation spaces. 

To be more precise, we mention here main results from \cite{OG}. Note that these results will be used to prove our new ones. 
We refer to Section \ref{section2} for definitions and notation.

The following theorem describes when the $K$-space $(X_0, X_1)_{0, q, b; K}$ with $b \in SV(0,\infty)$ 
coincides with the $J$-space $(X_0, X_1)_{0, q, a; J}$ 
with a  convenient $a \in SV(0,\infty).$

\begin{theorem}{\rm (cf. \cite [Theorem 3.1]{OG})} 
\label{equivalence theorem1}
Let $(X_0, X_1)$ be a compatible couple and $1\le q <\infty$. If $b \in SV(0,\infty)$ satisfies
\begin{equation}\label{prop_slow_var_funct_b1}
\int_x^\infty t^{-1}b^{\,q}(t)\, dt <\infty  \quad\text{for all \ }x>0, \qquad 
\int_0^\infty t^{-1}b^{\,q}(t)\, dt =\infty,
\end{equation}
and $a \in  SV(0,\infty)$ is such that 
\begin{equation}\label{def_slow_var_funct_a1}
a(x)=b^{-q/q\,'}(x)\int_x^\infty t^{-1}b^{\,q}(t)\, dt  \quad\text{for all\ }x>0,
\end{equation}
then
\begin{equation}\label{K_space=J_space1}
(X_0, X_1)_{0, q, b; K}=(X_0, X_1)_{0, q, a; J}.\text{\footnotemark$^)$}
\end{equation}
\end{theorem}
\footnotetext{$^{)}$ This means that $(X_0, X_1)_{0, q, b; K}=(X_0, X_1)_{0, q, a; J}$ and the norms of these spaces are equvivalent.}

Theorem \ref{equivalence theorem1} has been applied to establish the next assertion (cf. \cite [Theorem 3.2]{OG}).

\begin{theorem}[Density theorem]\label{density theorem1}
Let $(X_0, X_1)$ be a compatible couple and $1\le q<\infty$. If $b \in SV(0,\infty)$ satisfies
\eqref{prop_slow_var_funct_b1}, 
then  the space $X_0\cap X_1$ is dense in $(X_0, X_1)_{0, q, b; K}$.
\end{theorem}

We have also proved (see \cite [Theorems 3.3 and 3.4]{OG}) 
the following counterparts of the previous assertions.

\begin{theorem}
\label{equivalence theorem}
Let $(X_0, X_1)$ be a compatible couple and $1<q\le \infty$. If $a \in SV(0,\infty)$ satisfies
\begin{equation}\label{prop_slow_var_funct_a}
\int_0^x t^{-1}a^{-q'}(t)\, dt <\infty  \quad\text{for all \ }x>0, \qquad 
\int_0^\infty t^{-1}a^{-q'}(t)\, dt =\infty,
\end{equation}
and $b \in  SV(0,\infty)$ is defined by
\begin{equation}\label{def_slow_var_funct_b}
b(x):=a^{-q'/q}(x)\Big(\int_0^x t^{-1}a^{-q'}(t)\, dt \Big)^{-1} \quad\text{for all \ }x>0,
\end{equation}
then 
\begin{equation}\label{K_space=J_space}
(X_0, X_1)_{0, q, a; J}=(X_0, X_1)_{0, q, b; K}.
\end{equation}
\end{theorem}

\smallskip
\begin{theorem}[Density theorem]\label{density theorem}
Let $(X_0, X_1)$ be a compatible couple and $1<q<\infty$. If \,$a \in SV(0,\infty)$ satisfies
\eqref{prop_slow_var_funct_a}, 
then  the space $X_0\cap X_1$ is dense in $(X_0, X_1)_{0, q, a; J}$.
\end{theorem}

\smallskip
\begin{remark}\label{105} {\rm The first assumption in \eqref{prop_slow_var_funct_b1} guarantees that 
the space $(X_0, X_1)_{0, q, b; K}$ is an intermediate space between $X_0$ and $X_1$. Note also that 
if 
$\int_1^\infty t^{-1}b^{\,q}(t)\, dt =\infty$, 
then  $(X_0, X_1)_{0, q, b; K}=\{0\}$ (cf. Theorem \ref{302} mentioned below). On the other hand, if 
$\int_0^\infty t^{-1}b^{\,q}(t)\, dt <\infty$, then the space $(X_0, X_1)_{0, q, b; K}$ 
need not be described as a $(X_0, X_1)_{0, q, a; J}$ space with any $a \in SV(0,\infty)$ 
(cf. \cite [Lemma 4.8]{OG}).

Similarly, the first assumption in \eqref{prop_slow_var_funct_a} guarantees that 
the space $(X_0, X_1)_{0, q, b; J}$ is an intermediate space between $X_0$ and $X_1$. 
Note also that 
if 
$\int^1_0 t^{-1}a^{-q'}(t)\, dt =\infty$,  
then the functional 
$\|.\|_{\theta,q,v;J}$ vanishes on $X_0\cap X_1$ and thus it is not a norm provided that 
$X_0\cap X_1\ne \{0\}$ (cf. Theorem \ref{305} mentioned below).
}
\end{remark}

The next assertion corresponds to Theorems \ref{equivalence theorem1} 
with $q=\infty$.
\begin{theorem}[{\cite[Theorem 9.1]{OG}}]\label{equivalence theorem1*}
Let $(X_0, X_1)$ be a compatible couple. If $b \in SV(0,\infty)\cap AC(0,\infty)$ satisfies
\begin{equation}\label{prop_slow_var_funct_b1*}
b \text{\ is strictly decreasing},\quad b(0)=\infty, \quad b(\infty)=0,
\end{equation}
and if 
\begin{equation}\label{def_slow_var_funct_a1*}
a(x):= \frac{b^{\,2}(x)}{x (-b\, '(x))} \quad\text{for a.a. }x>0,
\end{equation}
then 
\begin{equation}\label{K_space=J_space1*}
(X_0, X_1)_{0, \infty, b; K}=(X_0, X_1)_{0, \infty, a; J}.
\end{equation}
\end{theorem}

\begin{remark}\label{1111*} {\rm Note that, by Lemma \ref{l2.6+} (iii) mentioned below, the function 
$b \in SV(0,\infty)$ is equivalent to $\overline{b} \in SV(0,\infty)\cap AC(0,\infty)$ given by $\overline{b}(t):= t^{-1}\int_0^t b(s)\,ds, \,t\in (0, \infty)$. One can also prove that  
if $b$ satisfies \eqref{prop_slow_var_funct_b1*}, then \eqref{prop_slow_var_funct_b1*} also holds with $b$ replaced by $\overline{b}$. 
Thus, Theorem \ref{equivalence theorem1*} remains true if the assumption $b \in SV(0,\infty)\cap AC(0,\infty)$ is replaced only by 
$b\in SV(0,\infty)$ provided that in \eqref{def_slow_var_funct_a1*} we write $\overline{b}$ instead of $b$, and $\overline{b}\,'$ instead of $b'$.}
\end{remark}

The next two theorems correspond to Theorems \ref{equivalence theorem} and~ 
\ref{density theorem} with $q=1$.
\begin{theorem}[{\cite[Theorem 10.1]{OG}}]\label{equivalence theorem1**}
Let $(X_0, X_1)$ be a compatible couple. If $a \in SV(0,\infty)\cap AC(0,\infty)$ satisfies
\begin{equation}\label{prop_slow_var_funct_a*}
a \text{\ is strictly decreasing},\quad a(0)=\infty, \quad a(\infty)=0,
\end{equation}
and if
\begin{equation}\label{def_slow_var_funct_b*}
b(x):= -x \,a\, '(x) \quad\text{for a.a. }x>0,
\end{equation}
then 
\begin{equation}\label{K_space=J_space1**}
(X_0, X_1)_{0, 1, a; J}=(X_0, X_1)_{0, 1, b; K}.
\end{equation}
\end{theorem}

\begin{remark}\label{1111} {\rm Similarly to Remark \ref{1111*}, one can show that there is a variant of Theorem~\ref{equivalence theorem1**}, where the assumption $a \in SV(0,\infty)\cap AC(0,\infty)$ is replaced only by $a \in SV(0,\infty)$. Details are left to the reader.
}
\end{remark}

\begin{theorem}[Density theorem, {\cite[Theorem 10.3]{OG}}]
\label{density theorem*}
Let $(X_0, X_1)$ be a compatible couple. If the function $a$ 
satisfies the assumptions of {\rm  Theorem \ref{equivalence theorem1**}}, then 
$X_0\cap X_1$ is dense in $(X_0, X_1)_{0, 1
, a; J}$.
\end{theorem}

We continue with our new results, which concern the cases when conditions in the above mentioned theorems 
are not satisfied but the spaces in question are still meaningful, i.e., they are Banach spaces and intermediate  spaces between $X_0$ and $X_1$. 

\smallskip
First we mention the next two assertions (which are  complements of Theorems \ref{equivalence theorem1} and \ref{density theorem1}), where we consider the case when condition \eqref{prop_slow_var_funct_b1} is not satisfied. 

\begin{theorem}
\label{KS}
Let $(X_0,X_1)$ be a compatible couple and $1\le q<\infty$. If $b\in SV\,(0,\infty)$ satisfies
\begin{equation}\label{KSB}
\int_0^\infty t^{-1}b^{\,q}(t)\, dt <\infty,  
\end{equation}
the function $B\in SV\,(0,\infty)$ is given by 
\begin{equation}\label{3131}
B(x):=b(x) \ \ \text{if} \ x\in [1,\infty),\qquad 
B(x)\approx \beta(x) \ \ \text{if} \ x\in (0,1),
\end{equation}
where $\beta\in \, SV (0,\infty)$ is such that
\begin{equation}\label{3141}
      \int_x^\infty t^{-1} \beta^q(t)\, dt <\infty \quad \text{for all} \ x>0, \quad \quad   
			\int_0^\infty t^{-1} \beta^q(t)\, dt=\infty, 
\end{equation}
and 
\begin{equation}\label{B1=}
A(x):=B^{-{q}/{q\,'}}(x)\int_x^{\infty}t^{-1}B^q(t)\,dt\quad \text{for all }x>0,\text{\footnotemark$^)$}
\end{equation}
then 
\begin{equation}\label{D2E1}
(X_0,X_1)_{0,q,b;K}=(X_0, X_0+X_1)_{0,q,B;K}=(X_0,X_0+X_1)_{0,q,A;J}.
\end{equation} 
\end{theorem}
\footnotetext{$^{)}$ If $q=1$, then $q'=\infty$ and so $A(x)=\int_x^{\infty}t^{-1}B(t)\,dt$ .}

\begin{theorem}[Density theorem]\label{density3}
Let $(X_0,X_1)$ be a compatible couple and $1\le q<\infty$. If $b\in SV\,(0,\infty)$ satisfies \eqref{KSB},
then the space $X_0 $ is dense in $(X_0, X_1)_{0, q, b; K}$.
\end{theorem}

\smallskip
Note that the space $X_0\cap X_1$ need not be dense in $(X_0, X_1)_{0, q, b; K}$ provided that \eqref{KSB}
holds (cf. \cite [Remark 3.9]{CS}).

\begin{remark}\label{remark100}\rm{Under the assumptions of Theorem \ref{KS},
\begin{equation}\label{101}
(X_0,X_1)_{0,q,b;K}=(X_0,X_1)_{0,q,b;K;(1,\infty)},
\end{equation}
\begin{equation}\label{1002}
(X_0,X_0+X_1)_{0,q,B;K}=(X_0,X_0+X_1)_{0,q,B;K;(1,\infty)},
\end{equation}
\begin{equation}\label{1003}
(X_0,X_0+X_1)_{0,q,A;J}=(X_0,X_0+X_1)_{0,q,A;J;(1,\infty)}.
\end{equation}

\noindent
Here, for a compatible couple $(Y_0, Y_1)$, $v \in SV(0,\infty)$, and $1\le q\le\infty$, we put
\begin{equation}\label{1004}
(Y_0,Y_1)_{0,q,v;K;(1,\infty)}:=\{f\in Y_0+Y_1: \|f\|_{0,q,v;K;(1, \infty)}<\infty\},
\end{equation}
where
\begin{equation}\label{1005}
 \|f\|_{0,q,v;K;(1,\infty)}:=\left\|t^{-1/q}\,v(t)\,K(f,t;Y_0,Y_1)\right\|_{q,(1,\infty)}.
\end{equation}

Similarly, the space $(Y_0,Y_1)_{0,q,v;J;(1,\infty)}$ consists of all $f \in Y_0+Y_1$ for which there is a~strongly measurable function $u :(1,\infty)\rightarrow Y_0 \cap Y_1$ such that 
\begin{equation}\label{1006}
f=\int_1^\infty u(s)\, \frac{ds}{s} \quad  {\rm (} {\text convergence\  in \ } Y_0+Y_1{\rm )}
\end{equation}
and for which the functional
\begin{equation}\label{1007}
\|f\|_{0,q,v;J;(1,\infty)}:=\inf \left\|t^{-1/q}\,v(t)\,J(u(t),t;Y_0,Y_1)\right\|_{q,(1,\infty)}
\end{equation} 
is finite {\rm (}the infimum extends over all representations \eqref{1006} of $f${\rm )}.  
} 
\end{remark}

\medskip
Note that \eqref{101} implies that \textit{the norm} of the space $(X_0,X_1)_{0,q,b;K;(1,\infty)}$ 
\textit{is equivalent to the norm} 
 of the space $(X_0,X_1)_{0,q,b;K}$. Similar assertions can be said about relations \eqref{1002} and \eqref{1003}. 
Analogous statements concern Remarks \ref{remark100J}, \ref{remark100*}, and \ref{remark100J1} mentioned below.

\smallskip
\begin{corollary}\label{cor2/2}
Under the assumptions of {\rm Theorem \ref{KS}},
\begin{equation}\label{corres2/2}
(X_0, X_1)_{0, q, b; K}=X_0 + (X_0,X_1)_{0, q, B; K}=X_0 + (X_0,X_1)_{0, q, A; J}.
\end{equation}
\end{corollary}

\medskip
In the following two assertions (which are complements of Theorems~\ref{equivalence theorem} and 
\ref{density theorem}) we consider the case when condition \eqref{prop_slow_var_funct_a} is not satisfied.

\begin{theorem}\label{JS11}
Let $(X_0, X_1)$ be a compatible couple and $1<q\le\infty$.
If $a \in SV(0,\infty)$ satisfies
\begin{equation}\label{akon111}
\int_0^{\infty} t^{-1}a^{-q\,'}(t)\, dt <\infty,
\end{equation}
the function $A\in SV(0,\infty)$ is given by
\begin{equation}\label{Afce111}
A(x):=a(x) \quad\text{if \ } x\in (0,1], \qquad 
A(x)\approx \alpha(x)\quad\text{if \ } x\in (1, \infty),
\end{equation}
where the function $\alpha \in SV(0,\infty)$ is such that
\begin{equation}\label{alfa111}
\int_0^x t^{-1}\alpha^{-q\,'}(t)\, dt <\infty\quad\text{for all \ }x>0, \qquad 
\int_0^\infty t^{-1}\alpha^{-q\,'}(t)\, dt =\infty,
\end{equation}
and
\begin{equation}\label{Bfce111}
B(x):=A^{-q\,'\!/q}(x) \left(\int_0^x t^{-1}A^{-q\,'}(t)\,dt\right)^{-1} \quad\text{for all \ }x>0,\text{\footnotemark$^)$}
\end{equation}
then
\begin{equation}\label{dual}
(X_0, X_1)_{0, q, a; J}=(X_0, X_0\cap X_1)_{0,q,A;J}=(X_0, X_0\cap X_1)_{0, q, B; K}.
\end{equation}
\end{theorem}
\footnotetext{$^{)}$ If $q=\infty$, then $q'=1$ and so $B(x)=\left(\int_0^x t^{-1}A^{-1}(t)\,dt\right)^{-1}$.}

\begin{theorem}[Density theorem]\label{density4}
Let $(X_0, X_1)$ be a compatible couple and $1<q<\infty$. If \,$a \in SV(0,\infty)$ satisfies
\eqref{akon111}, 
then  the space $X_0\cap X_1$ is dense in $(X_0, X_1)_{0, q, a; J}$.
\end{theorem}

\begin{remark}\label{remark100J}\rm{Note that, under the assumptions of Theorem \ref{JS11},
\begin{equation}\label{101J}
(X_0,X_1)_{0,q,a;J}=(X_0,X_1)_{0,q,a;J;(0,1)},
\end{equation}
\begin{equation}\label{1002J}
(X_0,X_0\cap X_1)_{0,q,A;J}=(X_0,X_0\cap X_1)_{0,q,A;J;(0,1)},
\end{equation}
\begin{equation}\label{1003J}
(X_0,X_0\cap X_1)_{0,q,B;K}=(X_0,X_0\cap X_1)_{0,q,B;K;(0,1)}.
\end{equation}

\noindent
Here, for a compatible couple $(Y_0, Y_1)$, $v \in SV(0,\infty)$  and $1\leq q \le\infty$, we put
\begin{equation}\label{1004J}
(Y_0,Y_1)_{0,q,v;K;(0,1)}:=\{f\in Y_0+Y_1: \|f\|_{0,q,v;K;(0,1)}<\infty\},
\end{equation}
where
\begin{equation}\label{1005J}
 \|f\|_{0,q,v;K;(0,1)}:=\left\|t^{-1/q}\,v(t)\,K(f,t;Y_0,Y_1)\right\|_{q,(0,1)}.
\end{equation}

Similarly, the space $(Y_0,Y_1)_{0,q,v;J;(0,1)}$ consists of all $f \in Y_0+Y_1$ for which there is a~strongly measurable function 
$u :(0,1)\rightarrow Y_0 \cap Y_1$ such that 
\begin{equation}\label{1006J}
f=\int_0^1 u(s)\, \frac{ds}{s} \quad  {\rm (} {\text convergence\  in \ } Y_0+Y_1{\rm )}
\end{equation}
and for which the functional
\begin{equation}\label{1007J}
\|f\|_{0,q,v;J;(0,1)}:=\inf \left\|t^{-1/q}\,v(t)\,J(u(t),t;Y_0,Y_1)\right\|_{q,(0,1)}
\end{equation} 
is finite {\rm (}the infimum extends over all representations \eqref{1006J} of $f${\rm )}.
} 
\end{remark}

\begin{corollary}\label{cor2}
Under the assumptions of {\rm Theorem \ref{JS11}},
\begin{equation}\label{corres2}
(X_0, X_1)_{0, q, a; J}=X_0 \cap (X_0,X_1)_{0, q, A; J}=X_0 \cap (X_0,X_1)_{0, q, B; K}.
\end{equation}
\end{corollary}

The next results is a complement of Theorem \ref{equivalence theorem1*}.

\begin{theorem}
\label{*equivalence theorem1*}
Let $(X_0, X_1)$ be a compatible couple and let $b \in SV(0,\infty)\cap AC(0,\infty)$ satisfy
\begin{equation}\label{*prop_slow_var_funct_b1*}
b \text{\ is strictly decreasing},\quad b(0)<\infty, \quad b(\infty)=0.
\end{equation}
Assume that the function $B\in SV(0,\infty)$ is given by 
\begin{equation}\label{3131*}
B(x):=b(x) \ \ \text{if} \ x\in [1,\infty),\qquad 
B(x):=c\,\beta(x) \ \ \text{if} \ x\in (0,1),
\end{equation}
where $\beta\in \, SV(0,\infty)\cap AC(0,\infty)$ is such that
\begin{equation}\label{aqinfty}
\beta \text{\ is strictly decreasing},\quad \beta(0)=\infty, \quad \beta(\infty)=0,      
\end{equation}
and $c$ is a positive constant chosen in such a way that $B \in AC(0,\infty)$. 
If 
\begin{equation}\label{difeq}
A(x):= \frac{B^{\,2}(x)}{x (-B\, '(x))} \quad\text{for a.a. }x>0,
\end{equation}
then
\begin{equation}\label{D2E1*}
(X_0,X_1)_{0,\infty,b;K}=(X_0, X_0+X_1)_{0,\infty,B;K}=(X_0,X_0+X_1)_{0,\infty,A;J}.
\end{equation} 
\end{theorem}

\begin{remark}\label{remark100*}\rm{Note that, under the assumptions of Theorem \ref{*equivalence theorem1*},
\begin{equation}\label{101*}
(X_0,X_1)_{0,\infty,b;K}=(X_0,X_1)_{0,\infty,b;K;(1,\infty)},
\end{equation}
\begin{equation}\label{1002*}
(X_0,X_0+X_1)_{0,\infty,B;K}=(X_0,X_0+X_1)_{0,\infty,B;K;(1,\infty)}.
\end{equation}
Moreover, if $\|A(t)\|_{\infty, (1,e)}<\infty$, then also
\begin{equation}\label{1003*}
(X_0,X_0+X_1)_{0,\infty,A;J}=(X_0,X_0+X_1)_{0,\infty,A;J;(1,\infty)}.
\end{equation}
}
\end{remark}

\smallskip
\begin{corollary}\label{cor2/3}
Under the assumptions of {\rm Theorem \ref{*equivalence theorem1*}},
\begin{equation}\label{corres2/3}
(X_0, X_1)_{0, \infty, b; K}=X_0 + (X_0,X_1)_{0, \infty, B; K}=X_0 + (X_0,X_1)_{0,\infty, A; J}.
\end{equation}
\end{corollary}

\smallskip
The following two theorems are complements of Theorems \ref{equivalence theorem1**} and  \ref{density theorem*}.     .

\begin{theorem}
\label{equivalence theorem1**a}
Let $(X_0, X_1)$ be a compatible couple and let $a \in SV(0,\infty)\cap AC(0,\infty)$ satisfy
\begin{equation}\label{prop_slow_var_funct_a**}
a \text{\ is strictly decreasing},\quad a(0)=\infty, \quad a(\infty)>0.
\end{equation}
Assume that the function $A\in SV(0,\infty)$ is given by 
\begin{equation}\label{3131*A}
A(x):=a(x) \ \ \text{if} \ x\in (0,1],\qquad 
A(x):=c\,\alpha(x) \ \ \text{if} \ x\in (1,\infty),
\end{equation}
where $\alpha\in \, SV(0,\infty)\cap AC(0,\infty)$ is such that
\begin{equation}\label{bq1}
\alpha \text{\ is strictly decreasing},\quad \alpha(0)=\infty, \quad \alpha(\infty)=0,      
\end{equation}
and $c$ is a positive constant chosen in such a way that $A \in AC(0,\infty)$. If 
\begin{equation}\label{difeqB}
B(x):=-x \,A\, '(x) \quad\text{for a.a. }x>0, 
\end{equation}
then
\begin{equation}\label{D2E1*S}
(X_0,X_1)_{0,1,a;J}=(X_0, X_0\cap X_1)_{0,1,A;J}=(X_0,X_0\cap X_1)_{0,1,B;K}.
\end{equation} 
\end{theorem}

\begin{theorem}[Density theorem]
\label{density theorem**}
Let $(X_0, X_1)$ be a compatible couple. If the function $a$ 
satisfies the assumptions of {\rm  Theorem \ref{equivalence theorem1**a}}, then 
$X_0\cap X_1$ is dense in $(X_0, X_1)_{0, 1, a; J}$.
\end{theorem}

\begin{remark}\label{remark100J1} \rm{Note that, under the assumptions of Theorem \ref{equivalence theorem1**a},
\begin{equation}\label{101J1}
(X_0,X_1)_{0,1,a;J}=(X_0,X_1)_{0,1,a;J;(0,1)},
\end{equation}
\begin{equation}\label{1002J1}
(X_0,X_0\cap X_1)_{0,1,A;J}=(X_0,X_0\cap X_1)_{0,1,A;J;(0,1)},
\end{equation}
\begin{equation}\label{1003J1}
(X_0,X_0\cap X_1)_{0,1,B;K}=(X_0,X_0\cap X_1)_{0,1,B;K;(0,1)}.
\end{equation}
} 
\end{remark}

\begin{corollary}\label{cor21} 
Under the assumptions of {\rm Theorem \ref{equivalence theorem1**a}},
\begin{equation}\label{corres21}
(X_0, X_1)_{0, 1, a; J}=X_0 \cap (X_0,X_1)_{0, 1, A; J}=X_0 \cap (X_0,X_1)_{0, 1, B; K}.
\end{equation}
\end{corollary}

\smallskip
Note that if the weight function $b$ is of logarithmic type, then the description of limiting K-spaces by means of limiting J-spaces
have been investigated in \cite{CK} (if the pair $(X_0,X_1)$ is ordered) and in \cite{CS}, \cite{BCFC20}
 (if $(X_0,X_1)$ is a general pair).

\smallskip
It makes also sense to consider variants of Theorem \ref{equivalence theorem1*}, where condition \eqref{prop_slow_var_funct_b1*}
is replaced by 
$$
b \text{\ is strictly decreasing},\quad b(0)=\infty, \quad b(\infty)>0,
$$
or by
$$
b \text{\ is strictly decreasing},\quad b(0)<\infty, \quad b(\infty)>0.
$$ 
(Recall that in Theorem \ref{*equivalence theorem1*} only the case 
$$
b \text{\ is strictly decreasing},\quad b(0)<\infty, \quad b(\infty)=0
$$
has been treated.) 
Then the corresponding space $(X_0, X_1)_{0, \infty, b; K}$ is also a Banach space and an intermediate space between $X_0$ and $X_1$ 
(cf. Theorem \ref{302} mentioned below). Similar remark can be made about Theorem \ref{equivalence theorem1**}. We postpone 
the investigation of these variants to another paper.

\smallskip
The paper is organized as follows: Section 2 contains notation, definitions and preliminaries.  
In Section~3 we collect auxiliary assertions. Sections 4--7 are devoted to proofs of our new main results. 
More precisely, the proofs of Theorems \ref{KS}, \ref{density3}, Remark \ref{remark100}, and Corollary \ref{cor2/2} 
are in Section 4, while Section 5 is devoted to the proofs of Theorems \ref{JS11}, \ref{density4}, Remark \ref{remark100J}, and Corollary \ref{cor2}. The proofs of Theorem \ref{*equivalence theorem1*}, Remark \ref{remark100*}, and Corollary \ref{cor2/3} 
can be find in Section 6. The last Section 7 contains proofs of Theorems \ref{equivalence theorem1**a}, 
\ref{density theorem**}, Remark \ref{remark100J1}, and Corollary \ref{cor21}.

\medskip
\section{Notation, definitions and preliminaries}\label{section2}

For two non-negative expressions ({i.e.} functions or functionals) ${\mcal
A}$, ${\mcal B}$, the symbol ${\mcal A}\lesssim {\mcal B}$ (or
${\mcal A}\gtrsim {\mcal B}$) means that $ {\mcal A}\leq c\, {\mcal
B}$ (or $c\,{\mcal A}\geq {\mcal B}$), where $c$ is a  positive constant independent of
significant quantities involved in $A$ and $B$. If ${\mcal A}\lesssim {\mcal
B}$ and ${\mcal A}\gtrsim{\mcal B}$, we write ${\mcal A}\approx
{\mcal B}$ and say that ${\mcal A}$ and ${\mcal B}$ are equivalent.
Throughout the paper we use the abbreviation $\text{LHS}(*)$
($\text{RHS}(*)$) for the left- (right-) hand side of relation
$(*)$. We adopt the convention that $a/(\infty)=0$, 
$a/0=\infty$ and $(\infty)^a=\infty$ for all $a\in(0,\infty)$. If $q\in[1,\infty]$, the conjugate
number $q\,'$ is given by $1/q+1/q\,'=1$. 
In the whole paper $\|.\|_{q;(c,d)},\, q\in
[1,\infty]$, denotes the usual $L_q$-norm on the
interval $(c,d)\subseteq \bkR$.

We denote by $\chi_{\Omega}$ the characteristic function of the set $\Omega$. 
The symbol ${\mcal M}^{+}(0,\infty)$ stands for the class of all (Lebesgue-) measurable  functions on the interval $(0, \infty)$, which are non-negative almost everywhere on $(0, \infty)$. We will admit only positive, finite weights, and thus we put 
$$
{\mcal W}(0,\infty):=\{w\in{\mcal M}^{+}(0,\infty): w<\infty \ \mbox{a.e. on}\  (0, \infty)\}.
$$
The symbol $AC(0,\infty)$ is used to denote the set of all absolutely continuous functions on the interval $(0,\infty)$.

If $f$ is a monotone function on the interval $(0,\infty)$, then we put 
$$
f(0):=\lim_{t\rightarrow 0_+} f(t) \qquad {\text and} \qquad
f(\infty):=\lim_{t\rightarrow \infty} f(t).
$$

We say that a positive, finite and
Lebesgue-measurable function $b$ is {\it slowly varying} on
$(0,\infty)$, and write $b\in SV(0,\infty)$, if, for each
$\varepsilon>0$, $t^{\varepsilon}b(t)$ is equivalent to a non-decreasing
function on $(0,\infty)$ and  $t^{-\varepsilon}b(t)$ is equivalent to
a~non-increasing function on $(0,\infty)$. Here we follow the definition of $SV(0,\infty)$
given in \cite{GOT2002Lrrinwsvf}; for other definitions see, for example,
\cite{BGT87:RV}, \cite{EEv04:HOFSE}, \cite{EKP:OISRIQ}, and \cite{Nev02:LKSBRPE}. The
family $SV(0,\infty)$ includes not only powers of iterated logarithms
and the broken logarithmic functions of \cite{EO00:RILFR} but also
such functions as $t\mapsto\exp\left(  \left\vert \log
t\right\vert ^{a}\right)  ,$ $a\in(0,1).$ (The last mentioned
function has the interesting property that it tends to infinity more
quickly than any positive power of the logarithmic function).

\smallskip

We mention some properties of slowly
varying functions. 
\begin{lemma} {\rm (\cite [Lemma 2.1]{OG})} \label{l2.6+}
Let $b,b_1,b_2 \in SV(0,\infty)$ and let $d(t):=b(1/t)$ for all $t>0 .$

\noindent
{\rm {\bf (i)}} Then $\, b_1 b_2, d \in SV(0,\infty)$ and $b^r \in SV(0,\infty)$ 
for each $r \in {\mathbb R}.
$ 

\noindent
{\rm {\bf (ii)}}  If $\,\varepsilon$ and $\kappa$ are positive numbers,
then
there are positive constants $c_\varepsilon$ and $C_\varepsilon $ \\
\hspace*{.8cm}such that
\[
c_\varepsilon \min \{ \kappa^{-\varepsilon}, \kappa^\varepsilon \}
b(t) \le b(\kappa t) \le C_\varepsilon \max \{ \kappa^\varepsilon,
\kappa^{-\varepsilon} \} b(t) \quad \mbox{for every }\, t>0 .
\]

\noindent
{\rm {\bf (iii)}} If $\,\alpha>0$ and $q \in (0,\infty]$, then, for all $t>0,$
\[
\|\tau^{\alpha-1/q} b(\tau)\|_{q,(0,t)} \approx t^\alpha b(t) \quad \mbox{and } \quad 
\|\tau^{-\alpha-1/q} b(\tau)\|_{q,(t,\infty)} \approx t^{-\alpha} b(t). 
\] 

\noindent
{\rm {\bf (iv)}} If $\,q \in (0,\infty]$, and     
$$  
B_0(t):=\|\tau^{-1/q} b(\tau)\|_{q,(0,t)}, \quad B_\infty(t):=\|\tau^{-1/q} b(\tau)\|_{q,(t,\infty)}, \ \ t>0,
$$
then
\begin{equation}\label{400}
b(t) \lesssim B_0(t),\quad b(t) \lesssim B_\infty(t), \ \ \text{for all \ } t>0. 
\end{equation}
Moreover, if $B_i(1) <\infty$, then $B_i$ 
belongs to $SV(0,\infty)$, $i=0, \infty$.

\noindent
{\rm {\bf (v)}} If $\,q \in (0,\infty)$, then
\begin{equation}\label{414}
\limsup_{t\rightarrow 0_+}\dfrac{\|\tau^{-1/q} b(\tau)\|_{q,(0,t)}}{b(t)}=\infty,
\quad  \quad
\limsup_{t\rightarrow \infty}\dfrac{\|\tau^{-1/q} b(\tau)\|_{q,(0,t)}}{b(t)}=\infty,
\end{equation}
\begin{equation}\label{415}
\limsup_{t\rightarrow 0_+}\dfrac{\|\tau^{-1/q} b(\tau)\|_{q,(t,\infty)}}{b(t)}=\infty,
\quad  \quad
\limsup_{t\rightarrow \infty}\dfrac{\|\tau^{-1/q} b(\tau)\|_{q,(t,\infty)}}{b(t)}=\infty.
\end{equation}

\end{lemma}


More properties and examples of slowly varying functions can be
found in \cite[Chapter~V, p. 186]{Zyg57:TS}, \cite{BGT87:RV}, 
\cite{VojislavMaric:RVDE00}, \cite{Nev02:LKSBRPE}, 
\cite{GOT2002Lrrinwsvf} and \cite{GNO10:PotAnal}.

\smallskip
Let $X$ and $Y$ be two Banach spaces. We say that $X$ 
{\it coincides} with $Y$ (and write $X=Y$) if $X$ and $Y$ are
equal in the algebraic and topological sense (their norms are equivalent). 
The symbol $X\hookrightarrow Y$
 means that $X\subset Y$ and the natural embedding of $X$ in $Y$ is continuous. 
The norm of this embedding is denoted by $\|Id\,\|_{X\rightarrow Y}$. 

A pair $(X_0,X_1)$ of Banach spaces $X_0$ and $X_1$ is called a {\it compatible couple} if there is a~Hausdorff topological vector space $\mathcal{X}$ in which each of $X_0$ and $X_1$ is continuously embedded.

If  $(X_0,X_1)$ is a compatible couple, then a Banach space $X$ is said to be an {\it intermediate space} 
between $X_0$ and $X_1$ if $X_0\cap X_1 \hookrightarrow X \hookrightarrow X_0+X_1.$

\begin{definition}\label{defin:1:interpolatonspaces}
Let $(X_0,X_1)$ be a compatible couple. 
\begin{trivlist}
\item[\hspace*{0.5cm}{\rm {\bf (i)}}] The Peetre $K$-functional is defined for each $f\in X_0+X_1$ and 
$t>0$ by
$$
K(f,t;X_0,X_1):=\inf\{\|f_0\|_{X_0}+t\|f_1\|_{X_1} : f=f_0+f_1\},
$$
where the infimum extends over all representations $f=f_0+f_1$ of $f$ with $f_0\in X_0$ and 
$f_1\in X_1$. 

\item[\hspace*{0.5cm}{\rm {\bf (ii)}}] The Peetre $J$-functional is defined for each $f\in X_0\cap X_1$ and 
$t>0$ by
$$
J(f,t;X_0,X_1):= \max \{ \left\|f\right\|_{X_0}, t\left\|f\right\|_{X_1} \}.
$$

\item[\hspace*{0.5cm}{\rm {\bf (iii)}}] For $0\leq \theta\leq 1$, $1\le q\le\infty$, and $v \in 
 {\mcal W}(0,\infty)$, 
we put
\begin{equation}\label{eq:circ_1}
(X_0,X_1)_{\theta,q,v;K}:=\{f\in X_0+X_1: \|f\|_{\theta,q,v;K}<\infty\},
\end{equation}
where
\begin{equation}\label{eq:circ_2}
 \|f\|_{\theta,q,v;K}\equiv\|f\|_{(X_0,X_1)_{\,\theta,q,v;K}}:=\left\|t^{-\theta-1/q}\,v(t)\,K(f,t;X_0,X_1)\right\|_{q,(0,\infty)}.
\end{equation}

\item[\hspace*{0.5cm}{\rm {\bf {(iv)}}}] Let $0\leq \theta\leq 1$, $1\le q\le\infty$, and let $v \in 
 {\mcal W}(0,\infty)$. The space $(X_0,X_1)_{\theta,q,v;J}$ consists of all $f \in X_0+X_1$ for which there is a strongly measurable function $u :(0,\infty)\rightarrow X_0 \cap X_1$ such that 
\begin{equation}\label{300}
f=\int_0^\infty u(s)\, \frac{ds}{s} \quad  {\rm (} {\text convergence\  in \ } X_0+X_1{\rm )}
\end{equation}
and for which the functional
\begin{equation}\label{301}
\|f\|_{\theta,q,v;J}\equiv\|f\|_{(X_0,X_1)_{\,\theta,q,v;J}}:=\inf \left\|t^{-\theta-1/q}\,v(t)\,J(u(t),t;X_0,X_1)\right\|_{q,(0,\infty)}
\end{equation} 
is finite {\rm (}the infimum extends over all representations \eqref{300} of $f${\rm )}.
\end{trivlist}
\end{definition}

We refer to Lemmas \ref{K} and \ref{J}  mentioned below 
for properties of the $K$-functional and the $J$-functional.

\begin{theorem}{\rm (\cite [Theorem 2.3]{OG})} \label{302}
Let $(X_0,X_1)$ be a compatible couple, $0\leq \theta\leq 1$, $1\le q\le\infty$, and let $v \in  {\mcal W}(0,\infty)$.

\noindent
{\rm {\bf A.}} If 
\begin{equation}\label{304}
\left\|t^{-\theta-1/q}\,v(t) \min\{1,t\}\right\|_{q,(0,\infty)}<\infty, 
\end{equation}
then:

{\rm {\bf (i)}} The space $(X_0,X_1)_{\theta,q,v;K}$ is an intermediate space between $X_0$ and $X_1$, that is,
$$
X_0 \cap X_1 \hookrightarrow (X_0,X_1)_{\theta,q,v;K} \hookrightarrow X_0+X_1.
$$

{\rm {\bf (ii)}} 
The space $ (X_0,X_1)_{\theta,q,v;K}$ is a Banach space.

\noindent
{\rm {\bf B.}}
If condition \eqref{304} is not satisfied, then $(X_0,X_1)_{\theta,q,v;K}=\{0\}.$
\end{theorem}


Note that assertion A of Theorem \ref {302} also follows from \cite [Proposition~3.3.1, p. 338]{BK}.

\begin{theorem}{\rm (\cite [Theorem 2.4]{OG})} \label{305}
Let $(X_0,X_1)$ be a compatible couple, $0\leq \theta\leq 1$, $1\le q\le\infty$, and let $v \in  {\mcal W}(0,\infty)$. 

\noindent
{\rm {\bf A.}} If  
\begin{equation}\label{307}
\left\|t^{\theta-1/q\,'}\,\frac{1}{v(t)} \min\left\{1,\frac{1}{t}\right\}\right\|_{q\,'\!,(0,\infty)}<\infty,
\end{equation}
then\,{\rm :}

{\rm {\bf (i)}} The space $(X_0,X_1)_{\theta,q,v;J}$ is an intermediate space between $X_0$ and $X_1$, that is,
$$
X_0 \cap X_1 \hookrightarrow (X_0,X_1)_{\theta,q,v;J} \hookrightarrow X_0+X_1.
$$

{\rm {\bf (ii)}} The space $ (X_0,X_1)_{\theta,q,v;J}$ is a Banach space.

\smallskip
\noindent
\noindent
{\rm {\bf B.}}
If condition \eqref{307} is not satisfied, 
then the functional 
$\|.\|_{\theta,q,v;J}$ vanishes on $X_0\cap X_1$ and thus it is not a norm provided that 
$X_0\cap X_1\ne \{0\}.$
\end{theorem}

\medskip
Recall also that the spaces $(X_0,X_1)_{0,q,v;K;(1,\infty)}, (X_0,X_1)_{0,q,v;J;(1,\infty)}$, 
$(X_0,X_1)_{0,q,v;K;(0,1)},$ $(X_0,X_1)_{0,q,v;J;(0,1)}$ have been defined in Remarks \ref{remark100} and \ref{remark100J} 
mentioned above.

\section{Auxiliary assertions}\label{auxi}

In the next  two lemmas some basic properties of the $K$- and $J$- functionals are summarized.

\begin{lemma} {\rm (cf. \cite [Proposition 1.2, p. 294]{BS:IO})} \label{about K}
If $(X_0,X_1)$ is a compatible couple, then, for each $f \in X_0+X_1$, the $K$-functional $K(f,t;X_0,X_1)$ is a nonnegative concave function of $t>0$, and 
\begin{equation}\label{K}
K(f,t;X_0,X_1)=t\,K(f,t^{-1};X_1,X_0) \quad \text{for all \ } t>0.
\end{equation}
In particular, $K(f,t;X_0,X_1)$ is non-decreasing on $(0,\infty)$ and $K(f,t;X_0,X_1)/t$ is non-increasing 
on $(0,\infty).$ 
\end{lemma}

\begin{lemma} {\rm (cf. \cite [Lemma 3.2.1, p. 42]{BL})}\label{about J}
If $(X_0,X_1)$ is a compatible couple, then, for each $f \in X_0\cap X_1$, the $J$-~functional $J(f,t;X_0,X_1)$ is a nonnegative convex function of $t>0$, and 
\begin{equation}\label{J}
J(f,t;X_0,X_1)=t\,J(f,t^{-1};X_1,X_0)\quad \text{for all \ } t>0,
\end{equation}
\begin{equation}\label{JJ}
J(f,t;X_0,X_1)\le \max\{1, t/s\}\,J(f,s;X_0,X_1) \quad \text{for all \ } t, s>0,
\end{equation}
\begin{equation}\label{KJ}
K(f,t;X_0,X_1)\le \min\{1, t/s\}\,J(f,s;X_0,X_1)\quad \text{for all \ } t, s>0.
\end{equation}
In particular, $J(f,t;X_0,X_1)$ is non-decreasing on $(0,\infty)$ and $J(f,t;X_0,X_1)/t$ is non-increasing 
on $(0,\infty).$ 
\end{lemma}

\smallskip
\begin{lemma}\label{vnoreni X_0 a X_1}
Let $(X_0, X_1)$ be a compatible couple.

\noindent
{\rm {\bf (i)}} If 
 $X_1\hookrightarrow X_0$ and $k:=\|Id\,\|_{X_1\rightarrow X_0}$, then 
\begin{equation}\label{K1}
K(f,t;X_0,X_1)=\|f\|_{X_0}\quad \text{for all} \ f \in X_0 +X_1=X_0 \ \text{and}\ \,t \ge k. 
\end{equation}

\noindent
{\rm {\bf (ii)}} If 
 $X_0\hookrightarrow X_1$ and and $k:=\|Id\,\|_{X_0\rightarrow X_1}$, then 
\begin{equation}\label{K0}
K(f,t;X_0,X_1)=t\,\|f\|_{X_1}\quad \text{for all} \ f \in X_0 +X_1=X_1 \ \text{and}\ \,t\in (0,1/k]. 
\end{equation}
\end{lemma}
\bpr
We refer to \cite[p. 47]{BL} for the proof of part (i). Part (ii) follows from part (i) and \eqref{K}.
\epr

\smallskip
\begin{corollary}\label{col}
{\rm {\bf (i)}} If all the assumptions of $\,${\rm Lemma \ref{vnoreni X_0 a X_1} (i)} are satisfied, then 
\begin{equation}\label{K1col}
K(f,t;X_0,X_1)\approx\|f\|_{X_0}\quad \text{for all} \ f \in X_0 +X_1=X_0 \ \text{and}\ \,t \ge 1. 
\end{equation}
{\rm {\bf (ii)}} If all the assumptions of $\,${\rm Lemma \ref{vnoreni X_0 a X_1} (ii)} are satisfied, then 
\begin{equation}\label{K0col}
K(f,t;X_0,X_1)\approx t\,\|f\|_{X_1}\quad \text{for all} \ f \in X_0 +X_1=X_1 \ \text{and}\ \,t\in (0,1). 
\end{equation}
\end{corollary}

\bpr
By Lemma \ref{vnoreni X_0 a X_1} (i), \eqref{K1col} holds if $k\le 1$.

Assume that $k>1$. 
Using \eqref{K1} and properties of the $K$-functional, we obtain, for all $f \in X_0$ and all $t \in [1, k)$, 
\begin{align*}
\|f\|_{X_0}=&K(f,k;X_0,X_1)\ge K(f,t;X_0,X_1) \ge t \frac{K(f,k;X_0,X_1)}{k}\\
 &\ge \frac{1}{k} K(f,k;X_0,X_1)\approx K(f,k;X_0,X_1)=\|f\|_{X_0},
\end{align*}
which, together with \eqref{K1}, implies \eqref{K1col}. 

Part (ii) is a consequence of part (i) and \eqref{K}.
\epr

\smallskip
\begin{lemma}\label{o K funct}
If $(X_0, X_1)$ is a compatible couple, then 
$$
K(f,t;X_0,X_0\cap X_1)\approx t\,\|f\|_{X_0}+K(f,t;X_0,X_1)\quad \text{for all}\ \,f\in X_0 \ \text{and}\ t\in (0,1). 
$$
\end{lemma}

\bpr
Let $f \in X_0, f=f_0+f_1,$ where $f_0 \in X_0$, $f_1 \in X_0\cap X_1$. Then 
\begin{equation}\label{K11}
K(f,t;X_0, X_0\cap X_1) \le \|f_0\|_{X_0} +t \,\|f_1\|_{X_0\cap X_1}\quad \text{for all }\ t \in (0, 1). 
\end{equation}
Moreover, 
\begin{equation}\label{K12}
{\rm RHS \eqref{K11}} \le \|f_0\|_{X_0} + t \left (\|f_1\|_{X_0} + \|f_1\|_{X_1}\right ).
\end{equation}
Since 
$$
\|f_1\|_{X_0}=\|f-f_0\|_{X_0}\le \|f\|_{X_0}+\|f_0\|_{X_0},
$$
we get, for all $t \in (0,1)$, 
\begin{align*}
{\rm RHS \eqref{K11}}&\le \|f_0\|_{X_0}+ t\,\left (\|f\|_{X_0}+\|f_0\|_{X_0} + \|f_1\|_{X_1}\right )\\
&\le 2 \left [t\,\|f\|_{X_0}+\left (\|f_0\|_{X_0}+t\, \|f_1\|_{X_1}\right )\right ]\\
&\approx t\,\|f\|_{X_0}+\left (\|f_0\|_{X_0}+t\, \|f_1\|_{X_1}\right ).
\end{align*}
Together with \eqref{K11}, this implies that
\begin{equation}\label{K13}
K(f,t;X_0, X_0\cap X_1)\lesssim  t\,\|f\|_{X_0}+K(f,t;X_0, X_1) \quad \text{for all}\ \,f\in X_0 \ \text{and}\ t\in (0,1).
\end{equation}

To prove the converse estimate for all $t\in (0,1)$, take $f \in X_0$ and assume that $f=f_0+f_1$, where $f_0 \in X_0, f_1 \in X_0 \cap X_1.$ 
Then, for all $t \in (0, 1)$,  
\begin{align*}
t\,\|f\|_{X_0}+K(f,t;X_0, X_1)&\le t\,\|f\|_{X_0}+\left (\|f_0\|_{X_0} + t\,\|f_1\|_{X_1}\right )\\
&\le t \left (\|f_0\|_{X_0} + \|f_1\|_{X_0}\right )+\left (\|f_0\|_{X_0} +t\, \|f_1\|_{X_1}\right )\\
&\le 2\, \left ( \|f_0\|_{X_0}+t \left(\|f_1\|_{X_0}+\|f_1\|_{X_1}\right )\right )\\
&\approx  \|f_0\|_{X_0}+t\, \|f_1\|_{X_0\cap X_1},
\end{align*}
which implies that, for all $f\in X_0 $ and all $t\in (0,1)$,
\begin{equation}\label{K14}
t\,\|f\|_{X_0}+K(f,t;X_0, X_1)\lesssim K(f,t;X_0, X_0\cap X_1).
\end{equation}
\epr

\begin{lemma}\label{1 o K funct}
If $(X_0, X_1)$ is a compatible couple, then 
$$
K(f,t;X_0,X_0 + X_1)= \min\{1,t\}\,K(f,\max\{1,t\};X_0,X_1)\quad \text{for all}\ \,f\in X_0+X_1 \ \text{and}\ t>0.
$$
\end{lemma}
\bpr The result follows from \cite[Thm. 2, p. 178]{Mal84} (see also \cite[proof of Lemma~4.6]{OG}.)
\epr

\begin{lemma}\label{1 o J funct}
If $(X_0, X_1)$ is a compatible couple, then 
$$
J(f,t;X_0,X_0 \cap X_1)= \max\{1,t\}\,J(f,\min\{1,t\};X_0,X_1)\quad \text{for all}\ \,f\in X_0\cap X_1 \ \text{and}\ t>0.
$$
\end{lemma}
\bpr Let $f \in X_0\cap X_1$. If $t>0$, then 
$$
J(f,t;X_0,X_0 \cap X_1)=\max\{\|f\|_{X_0}, t\|f\|_{X_0\cap X_1},\}
=\max\{\|f\|_{X_0}, t\|f\|_{X_0}, t\|f\|_{X_1}\}.
$$
Hence, 
$$
J(f,t;X_0,X_0 \cap X_1)=\max\{\|f\|_{X_0}, t\,\|f\|_{X_1}\}= J(f,t;X_0,X_1)\ \ \text{for all} \ t\in (0,1),
$$
$$
J(f,t;X_0,X_0 \cap X_1)=\max\{t\,\|f\|_{X_0}, t\,\|f\|_{X_1}\}
= t\,J(f,1;X_0,X_1)\ \ \text{for all} \ t\ge 1,
$$
and the result follows.
\epr

\begin{theorem}{\rm (\cite [Theorem 4.5]{OG})}\label{308}
Let $(X_0,X_1)$ be a compatible couple and $1\le q\le \infty.$ If 
$v \in  {\mcal W}(0,\infty)$ is such that 
$$
\left\| t^{-1/q}v(t)\right\|_{q, (0,\infty)}<\infty,
$$
then 
$$
\|f\|_{(X_0,X_1)_{0,q,v;K}}\approx \left\|t^{-1/q}\,v(t)\,K(f,t;X_0,X_1)\right\|_{q,(1,\infty)}
\quad \text {for\  all \ } f \in X_0+X_1.
$$
\end{theorem}

\smallskip
\begin{theorem}{\rm (\cite [Theorem 4.6]{OG})}\label{311} If all the assumptions of \,{\rm Lemma \ref{308}} are satisfied, then
$$
(X_0,X_1)_{0,q,v;K}=(X_0,X_0+X_1)_{0,q,v;K}.
$$
\end{theorem}

\smallskip
\begin{lemma}{\rm (\cite [Lemma 11.2]{OG})}\label{102} Assume that $1<q<\infty$.

\noindent
\!{\bf \,(i)} Let $b \in SV(0,\infty)$ be such that
\begin{equation}\label{prop_slow_var_funct_b1r}
\int_x^\infty t^{-1}b^{\,q}(t)\, dt <\infty  \quad\text{for all \ }x>0 \qquad \text{and\ \ }
\int_0^\infty t^{-1}b^{\,q}(t)\, dt =\infty.
\end{equation}
If 
\begin{equation}\label{def_slow_var_funct_a1r}
a(x):=b^{-q/q\,'}(x)\int_x^\infty t^{-1}b^{\,q}(t)\, dt  \quad\text{for all \ }x>0,
\end{equation}
then $a \in SV(0,\infty)$, 
\begin{equation}\label{prop_slow_var_funct_ar}
\int_0^x t^{-1}a^{-q\,'}(t)\, dt <\infty  \quad\text{for all \ }x>0 \qquad \text{and\ \ }
\int_0^\infty t^{-1}a^{-q\,'}(t)\, dt =\infty.
\end{equation}
Moreover,
\begin{equation}\label{def_slow_var_funct_br}
b(x)\approx a^{-q\,'\!/q}(x)\Big(\int_0^x t^{-1}a^{-q\,'}(t)\, dt \Big)^{-1} \quad\text{for a.a. \ }x>0
\end{equation}
and
\begin{equation}\label{103}
\Big(\int_0^x t^{-1}a^{-q\,'}(t)\, dt \Big)^{1/q\,'} \Big(\int_x^\infty t^{-1}b^{\,q}(t)\, dt \Big)^{1/q}
=\Big(\frac{1}{q\,'-1}\Big)^{1/q\,'} \quad\text{for all \ }x>0.
\end{equation}
\!{\bf \,(ii)} Let $a \in SV(0,\infty)$ be such that \eqref{prop_slow_var_funct_ar} holds. If 
\begin{equation}\label{def_slow_var_funct_br=}
b(x):=a^{-q\,'\!/q}(x)\Big(\int_0^x t^{-1}a^{-q\,'}(t)\, dt \Big)^{-1} \quad\text{for all \ }x>0,
\end{equation}
then $b \in SV(0,\infty)$ and \eqref{prop_slow_var_funct_b1r} is satisfied. Moreover, 
\begin{equation}\label{def_slow_var_funct_a1r=}
a(x)\approx b^{-q/q\,'}(x)\int_x^\infty t^{-1}b^{\,q}(t)\, dt  \quad\text{for a.a. \ }x>0
\end{equation}
and 
\begin{equation}\label{103*}
\Big(\int_0^x t^{-1}a^{-q\,'}(t)\, dt \Big)^{1/q\,'} \Big(\int_x^\infty t^{-1}b^{\,q}(t)\, dt \Big)^{1/q}
=\Big(\frac{1}{q-1}\Big)^{1/q} \quad\text{for all \ }x>0.
\end{equation}
\end{lemma}

\medskip
\section{Proofs of Theorems \ref{KS}, \ref{density3}, Remark \ref{remark100}, and Corollary \ref{cor2/2}}\label{Pr1Main}

{\bf I.} Proof of Theorem \ref{KS}. 
By Theorem \ref{308},
\begin{equation}\label{310}
\|f\|_{(X_0,X_1)_{0,q,b;K}}\approx \left\|t^{-1/q}\,b(t)\,K(f,t;X_0,X_1)\right\|_{q,(1,\infty)}
\quad \text {for\  all \ } f \in X_0+X_1,
\end{equation}
and
\begin{equation}\label{316}
\|f\|_{(X_0,X_0+X_1)_{0,q,b;K}}\approx 
\left\|t^{-1/q}\,b(t)\,K(f,t;X_0,X_0+X_1)\right\|_{q,(1,\infty)}
\quad \text {for\  all \ } f \in X_0+X_1.
\end{equation}
Moreover, by Theorem \ref{311},
\begin{equation}\label{312}
(X_0,X_1)_{0,q,b;K}=(X_0,X_0+X_1)_{0,q,b;K}.
\end{equation}
Using \eqref{KSB}-\eqref{3141} and properties of slowly varying functions, 
one can see that the function $B$ is satisfies
$$
\left\|t^{-1/q}\,B(t) \min\{1,t\}\right\|_{q,(0,\infty)}<\infty. 
$$
Thus, the space $(X_0, X_0+X_1)_{0,q,B;K}$ is a Banach space, which is intermediate between $X_0$ and $X_0+X_1$ (cf. Theorem \ref{302}).
Similarly, the space  $(X_0, X_0+X_1)_{0,q,b;K}$ is intermediate between $X_0$ and $X_0+X_1$, and hence
\begin{equation}\label{3031}
(X_0,X_0+X_1)_{0,q,b;K} \hookrightarrow X_0+X_1.
\end{equation}
As $X_0 \hookrightarrow X_0+X_1$, Corollary \ref{col} (ii) (with $X_0, X_1$ replaced by $X_0, X_0+X_1$) implies that
\begin{equation}\label{K0col1}
K(f,t;X_0,X_0+X_1)\approx t\,\|f\|_{X_0+X_1}\quad \text{for all} \ f \in X_0 +X_1 \ \text{and}\ \,t\in (0,1). 
\end{equation}
Using \eqref{3131}, \eqref{K0col1}, properties of slowly varying functions, \eqref{316}, and \eqref{3031}, we obtain 
\begin{align*}
&\|f\|_{(X_0,X_0+X_1)_{0,q,B;K}}\\
\quad &\approx 
\left\|t^{-1/q}\,\beta(t)\,K(f,t;X_0,X_0+X_1)\right\|_{q,(0,1)} + \left\|t^{-1/q}\,b(t)\,K(f,t;X_0,X_0+X_1)\right\|_{q,(1,\infty)}\\
\quad &\approx \left\|t^{1-1/q}\,\beta(t)\right\|_{q,(0,1)} \|f\|_{X_0+X_1}+\left\|t^{-1/q}\,b(t)\,K(f,t;X_0,X_0+X_1)\right\|_{q,(1,\infty)}\\
\quad&\approx \|f\|_{X_0+X_1}+\|f\|_{(X_0,X_0+X_1)_{0,q,b;K}}\\
\quad&\approx \|f\|_{(X_0,X_0+X_1)_{0,q,b;K}}\qquad \text{for all } \  f \in X_0 +X_1.
\end{align*}
Consequently, 
\begin{equation}\label{3121}
(X_0, X_0+X_1)_{0,q,b;K}=(X_0, X_0+X_1)_{0,q,B;K}.
\end{equation}
Making use of \eqref{312} and \eqref{3121}, we arrive at 
\begin{equation}\label{31211}
(X_0, X_1)_{0,q,b;K}=(X_0, X_0+X_1)_{0,q,B;K}.
\end{equation}
Moreover, the function $B$ satisfies 
\begin{equation}\label{B*}
\int_x^\infty t^{-1}B^q(t)\, dt <\infty \quad\text{for all }x>0,\qquad \int_0^\infty t^{-1}B^q(t)\, dt =\infty.
\end{equation}
Thus, applying Theorem \ref{equivalence theorem1} (with $(X_0, X_1), b,$ and $a$ replaced by $(X_0, X_0+X_1), B,$  and $A$, 
respectively), we arrive at 
\begin{equation}\label{B**}
(X_0, X_0+X_1)_{0,q,B;K}=(X_0, X_0+X_1)_{0,q,A;J}.
\end{equation}
Together with \eqref{31211}, this gives \eqref{D2E1}.
\hskip 8,6cm$\square$

\smallskip
{\bf II.} Proof of Theorem \ref{density3}. The function $B$ defined in Theorem \ref{KS} satisfies \eqref{B*}. Thus, applying 
Theorem \ref{density theorem1} (with $(X_0, X_1)$ replaced by $(X_0, X_0 + X_1)$ and with $b$ replaced by $B$) we obtain that 
the space $X_0=X_0 \cap (X_0 + X_1)$ is dense in $(X_0, X_0+X_1)_{0,q,B;K}$. Together with \eqref{D2E1}, this yields the result. 
\hskip 9cm$\square$

\medskip
{\bf III.} Proof of Remark \ref{remark100}. 
Relation \eqref{101} was already verified in the proof of Theorem~\ref{KS} (cf. \eqref{310}). 
Moreover, \eqref{1002} is a consequence of \eqref{3121}, \eqref{316}, and \eqref{3131}. 
Thus, it remains to prove that
\begin{equation}\label{10003}
Y:=(X_0,X_0+X_1)_{0,q,A;J}=(X_0,X_0+X_1)_{0,q,A;J;(1,\infty)}=:Y_{(1, \infty)}.
\end{equation} 

{\bf (i)} Assume that $f\in Y$. 
Then there is a strongly measurable function 
$u :(0,\infty)\rightarrow X_0=X_0 \cap (X_0+X_1)$ such that 
\begin{equation}\label{3000}
f=\int_0^\infty u(t)\, \frac{dt}{t} \quad  {\rm (} {\text convergence\  in \ } X_0+X_1{\rm )}
\end{equation}
and 
\begin{equation}\label{3001}
\left\|t^{-1/q}\,A(t)\,J(u(t),t;X_0,X_0+X_1)\right\|_{q,(0,\infty)}\le 2\, \left\|f\right\|_Y.
\end{equation}
Moreover,
\begin{equation}\label{3002}
f=\int_0^1 u(t)\, \frac{dt}{t}+\int_1^\infty u(t)\, \frac{dt}{t}=:f_0+f_\infty.
\end{equation}
First we show that $f_0\in X_0$. By Lemma \ref{102} (i) (with $b,a$ replaced by $B,A$), 
if $1<q<\infty$, and by \eqref{B*}, if $q=1$,
\text{\footnotemark$^)$}
\footnotetext{$^{)}$ Note that if $q=1$, then, by \eqref{B1=}, $A(x)=\int_x^{\infty}t^{-1}B(t)\,dt, x>0$.}
\begin{equation}\label{3003}
\|t^{-1/q'}A^{-1}(t)\|_{q',(0,1)} <\infty \quad \text{and}  \quad \|t^{-1/q'}A^{-1}(t)\|_{q',(0,\infty)} =\infty.
\end{equation}
Using \eqref{3002}, the definition of the J-functional, H\" older's inequality, \eqref{3003}, and \eqref{3001}, we obtain 
\begin{align*}
\|f_0\|_{X_0}&\le \int_0^1\|u(t)\|_{X_0}\,\frac{dt}{t}\le \int_0^1 J(u(t),t;X_0, X_0+X_1)\,\frac{dt}{t}\\
&\le \left\|t^{-1/q}\,A(t)\,J(u(t),t;X_0,X_0+X_1)\right\|_{q,(0,1)}\|t^{-1/q'}A^{-1}(t)\|_{q',(0,1)}\\
&\approx \left\|t^{-1/q}\,A(t)\,J(u(t),t;X_0,X_0+X_1)\right\|_{q,(0,1)}\le 2\,\|f\|_Y <\infty.
\end{align*}
Hence,
\begin{equation}\label{3004}
f_0\in X_0 \qquad \text{and}  \qquad \|f_0\|_{X_0}\lesssim  \|f\|_Y.
\end{equation}

Now we prove that 
\begin{equation}\label{3005}
X_0\hookrightarrow Y_{(1, \infty)}.
\end{equation}
If $g\in X_0$, then
$$
g=\int_1^\infty g\chi_{(1,e)}(t)\,\frac{dt}{t},\qquad g\chi_{(1,e)}:(1,\infty)\rightarrow X_0,
$$
which, together with the definition of the J-functional, 
the embedding $X_0 \hookrightarrow X_0+X_1$, and properties of slowly varying functions, implies that
$$
J(g\chi_{(1,e)}(t),t; X_0, X_0+X_1) \lesssim \|g\|_{X_0}\, \chi_{(1,e)}(t)\quad \text{for all } t>1. 
$$
Consequently,
\begin{align*}
\|g\|_{Y_{(1, \infty)}} &\le \left\|t^{-1/q}\,A(t)\,J(g\chi_{(1,e)}(t),t; X_0, X_0+X_1)\right\|_{q,(1, \infty)}\\
 &\lesssim \|g\|_{X_0}\,\|t^{-1/q}A(t)\|_{q,(1,e)}\approx \|g\|_{X_0},
\end{align*}
and \eqref{3005} is verified.

Combining \eqref{3004} and \eqref{3005}, we arrive at
\begin{equation}\label{3006}
f_0  \in Y_{(1,\infty)} \qquad \text{and} \qquad \|f_0\|_{Y_{(1,\infty)}} \lesssim \|f\|_Y.
\end{equation}

Since 
$$
f_\infty=\int_1^\infty u(t)\, \frac{dt}{t}\quad \text{and}\quad   u :(1,\infty)\rightarrow X_0,
$$
and, by \eqref{3001}, 
$$
\left\|t^{-1/q}\,A(t)\,J(u(t),t;X_0,X_0+X_1)\right\|_{q,(1,\infty)}\le 2\,\|f\|_Y,
$$
we get from the definition of the norm in the space $Y_{(1,\infty)}$ that 
$$
f_\infty  \in Y_{(1,\infty)} \qquad \text{and} \qquad \|f_\infty\|_{Y_{(1,\infty)}} \lesssim \|f\|_Y.
$$
Together with \eqref{3006} and the fact that $f=f_0+f_\infty$, this shows that 
$$
f \in Y_{(1,\infty)}   \qquad \text{and} \qquad \|f\|_{Y_{(1,\infty)}} \lesssim \|f\|_Y.
$$
Consequently, 
\begin{equation}\label{3007*}
Y\hookrightarrow Y_{(1, \infty)}.
\end{equation}

{\bf (ii)} Assume now that $f\in Y_{(1, \infty)}$. 
Then there exists a strongly measurable function 
$u:(1,\infty)\rightarrow X_0$ such that 
$$
f=\int_1^\infty u(t)\, \frac{dt}{t} \quad  {\rm (} {\text convergence\  in \ } X_0+X_1{\rm )}
$$
and 
$$
\left\|t^{-1/q}\,A(t)\,J(u(t),t;X_0,X_0+X_1)\right\|_{q,(1,\infty)}\le 2\, \left\|f\right\|_{Y_{(1,\infty)}}.
$$
On putting 
$$
w(t):=\begin{cases}0, &t \in (0, 1]\\
             u(t) \quad &t \in (1, \infty), 
		\end{cases}
$$
we see that 
$$
f=\int_0^\infty w(t)\, \frac{dt}{t} \quad  {\rm (} {\text convergence\  in \ } X_0+X_1{\rm )}
$$
and 
\begin{align*}
\left\|t^{-1/q}\,A(t)\,J(w(t),t;X_0,X_0+X_1)\right\|_{q,(0,\infty)}
&=\left\|t^{-1/q}\,A(t)\,J(u(t),t;X_0,X_0+X_1)\right\|_{q,(1,\infty)}\\
&\le  2\, \left\|f\right\|_{Y_{(1,\infty)}},
\end{align*}
which (together with the definition of the norm in the space $Y$) implies that 
\begin{equation}\label{3007**}
Y_{(1,\infty)}\hookrightarrow Y.
\end{equation} 

Finally, combining embeddings \eqref{3007*} and \eqref{3007**}, we get \eqref{10003}.
\hskip 3,92cm$\square$

\medskip
{\bf IV.} Proof of Corollary \ref{cor2/2}. 
By \eqref{D2E1},
\begin{equation}\label{D2E1/2}
(X_0,X_1)_{0,q,b;K}=(X_0, X_0+X_1)_{0,q,B;K}.
\end{equation} 
The fact that the space $(X_0, X_0+X_1)_{0,q,B;K}$ is an intermediate space between $X_0$ and $X_0+X_1$ 
means that
\begin{equation}\label{intermed}
 X_0=X_0\cap (X_0+X_1) \hookrightarrow (X_0, X_0+X_1)_{0,q,B;K} \hookrightarrow X_0+X_1.
\end{equation}
Since 
$$
K(f,t;X_0,X_1)\ge K(f,t;X_0,X_0+X_1) \quad \text{for all } f\in X_0+X_1 \ \text{and } t>0,
$$
we see that 
\begin{equation}\label{emb1}
(X_0,X_1)_{0,q,B;K}  \hookrightarrow (X_0, X_0+X_1)_{0,q,B;K}.
\end{equation}
Using the first embedding in \eqref{intermed} and \eqref{emb1},
we arrive at 
\begin{equation}\label{leftemb}
X_0+(X_0,X_1)_{0,q,B;K}  \hookrightarrow (X_0, X_0+X_1)_{0,q,B;K}.
\end{equation}

Now we are going to prove the opposite embedding to \eqref{leftemb}. Let $f \in X_0+X_1$ 
and let $f=f_0+f_1$ be such a decomposition of $f$ that 
$$
\|f_0\|_{X_0} + \|f_1\|_{X_1} \le 2\,\|f\|_{X_0+X_1}.
$$
Together with the second embedding in \eqref{intermed}, this shows that 

\begin{equation}\label{emb2}
\|f_0\|_{X_0} + \|f_1\|_{X_1} \le 2\,\|f\|_{(X_0, X_0+X_1)_{0,q,B;K}}.
\end{equation}

{\bf (i)} Assume that $t>1$. Using properties of the $K$-functional, Lemma \ref{1 o K funct}, and 
\eqref{emb2}, we obtain
\begin{align*}
K(f_1,t;X_0,X_1)&\le K(f,t;X_0,X_1)+K(f_0,t;X_0,X_1)\\
&=K(f,t;X_0,X_0+X_1)+K(f_0,t;X_0,X_1)\\
&\le K(f,t;X_0,X_0+X_1)+\|f_0\|_{X_0}\\
&\lesssim K(f,t;X_0,X_0+X_1)+\|f\|_{(X_0, X_0+X_1)_{0,q,B;K}}.
\end{align*}
Consequently, 
\begin{align*}
\|t^{-1/q}&B(t)K(f_1,t;X_0,X_1)\|_{q,(1,\infty)}\\
&\lesssim \|t^{-1/q}B(t)K(f,t;X_0,X_0+X_1)\|_{q,(1,\infty)}
+\|t^{-1/q}B(t)\|_{q,(1,\infty)} \|f\|_{(X_0, X_0+X_1)_{0,q,B;K}}.
\end{align*}
Since, by \eqref{3131} and \eqref{KSB},
$$
\|t^{-1/q}B(t)\|_{q,(1,\infty)}=\|t^{-1/q}b(t)\|_{q,(1,\infty)}<\infty,
$$
we obtain that
\begin{equation}\label{jednaner}
\|t^{-1/q}B(t)K(f_1,t;X_0,X_1)\|_{q,(1,\infty)}\lesssim \|f\|_{(X_0, X_0+X_1)_{0,q,B;K}}.
\end{equation}

{\bf (ii)} Let $t \in (0,1)$. Using the inequality
$$
K(f_1,t;X_0,X_1)\le t\,\|f_1\|_{X_1},
$$
properties of slowly varying functions, and \eqref{emb2}, we get
\begin{align}\label{druhaner}
\|t^{-1/q}B(t)K(f_1,t;X_0,X_1)\|_{q,(0,1)}&\le \|t^{1-1/q}B(t)\|_{q,(0,1)} \|f_1\|_{X_1}\\
& \lesssim \|f\|_{(X_0, X_0+X_1)_{0,q,B;K}}.\notag
\end{align}
Estimates \eqref{jednaner} and \eqref{druhaner} imply that 
$$
\|f_1\|_{(X_0,X_1)_{0,q,B;K}}\lesssim 
\|f\|_{(X_0, X_0+X_1)_{0,q,B;K}}.
$$
Together with \eqref{emb2}, this gives
$$
\|f_0\|_{X_0} + \|f_1\|_{(X_0,X_1)_{0,q,B;K}}\lesssim \|f\|_{(X_0, X_0+X_1)_{0,q,B;K}}
$$
for all $f \in X_0+X_1$. Consequently,
$$
(X_0, X_0+X_1)_{0,q,B;K}  \hookrightarrow  X_0+(X_0,X_1)_{0,q,B;K}.
$$
This embedding and \eqref{leftemb}  show that
\begin{equation}\label{rovnost}
(X_0, X_0+X_1)_{0,q,B;K} =  X_0+(X_0,X_1)_{0,q,B;K}.
\end{equation}
Making use of \eqref{D2E1/2} and \eqref{rovnost}, we arrive at
\begin{equation}\label{rovnost2}
(X_0,X_1)_{0,q,b;K}=X_0+(X_0,X_1)_{0,q,B;K}.
\end{equation}
Now, since \eqref{B*} holds, we can apply Theorem \ref{equivalence theorem1} (with $B$ and $A$ instead of $b$ and $a$, respectively) 
to get 
$$
(X_0,X_1)_{0,q,B;K}=(X_0,X_1)_{0,q,A;J},
$$
which, together with \eqref{rovnost2}, gives \eqref{corres2/2}.
\hskip 8cm$\square$

\medskip
\section{Proofs of Theorems \ref{JS11}, \ref{density4}, Remark \ref{remark100J}, and Corollary \ref{cor2}}\label{Pr1MainJ}

{\bf I.} Proof of Theorem \ref{JS11}. 
Put 
$$
Y:=(X_0, X_1)_{0, q, a; J}\qquad \text{and}\qquad Z:=(X_0, X_0\cap X_1)_{0, q, A; J}.
$$

{\bf (i)} Let $f \in Y$. Then there is a function $u:(0, \infty)\rightarrow X_0 \cap X_1$ such that
\begin{equation}\label{f_funkce}
f=\int_0^\infty u(t) \frac{dt}{t}   \quad  {\rm (} {\text convergence\  in \ } X_0+X_1{\rm )}
\end{equation}
and
\begin{equation}\label{nernorm}
\|t^{-1/q}a(t)J(u(t),t;X_0,X_1)\|_{q,(0,\infty)}\le 2\,\|f\|_Y.
\end{equation}

Set 
\begin{equation}\label{f_infty_funkce}
f_{\infty}:=\int_1^\infty u(t) \frac{dt}{t}, 
\end{equation}
\begin{equation}\label{w}
w(t):=\begin{cases} u(t)+t\,f_\infty, \quad &t \in (0, 1)\\
             0, &t\ge 1.
\end{cases}
\end{equation}
Then
\begin{equation}\label{w1}
f=\int_0^\infty w(t) \frac{dt}{t}\, .
\end{equation}
Moreover, we claim that
\begin{equation}\label{w v  pruniku}
w(t) \in  X_0 \cap X_1\quad \text{for all\ } t \in (0,\infty).
\end{equation}
 With respect to \eqref{w}, to verify \eqref{w v  pruniku}, it is sufficient to show that 
\begin{equation}\label{f_infty_funkce1}
f_{\infty} \in X_0 \cap X_1.
\end{equation}
Using \eqref{f_infty_funkce}, we get
$$
\|f_\infty\|_{X_0 \cap X_1} \le \int_1^\infty \|u(t)\|_{X_0 \cap X_1}\frac{dt}{t}
 \approx \int_1^\infty \left(\|u(t)\|_{X_0} + \|u(t)\|_{X_1}\right) \frac{dt}{t}.
$$
Since, for all $t>0$,  
$$
\|u(t)\|_{X_0} \le J(u(t),t;X_0, X_1)\qquad \text{and}\qquad \|u(t)\|_{X_1} \le \frac{J(u(t),t;X_0, X_1)}{t},
$$
we see that
$$
\|f_\infty\|_{X_0 \cap X_1} \lesssim \int_1^\infty J(u(t),t;X_0, X_1)\frac{dt}{t} +\int_1^\infty J(u(t),t;X_0, X_1)\frac{dt}{t^2}.
$$
Thus, applying H\" older's inequality, the fact that $\|t^{-1/q\,'}a^{-1}(t)\|_{q\,'\!,\, (1,\infty)} < \infty$ (cf. \eqref{akon111}), properties of slowly varying functions, and \eqref{nernorm}, we obtain
\begin{align*}
\|f_\infty\|_{X_0 \cap X_1} &\lesssim \|t^{-1/q}a(t)J(u(t),t;X_0,X_1)\|_{q,(1,\infty)} \,\|t^{-1/q\,'}a^{-1}(t)\|_{q\,'\!,\, (1,\infty)}\\
&\ \ +\|t^{-1/q}a(t)J(u(t),t;X_0,X_1)\|_{q,(1,\infty)} \,\|t^{-1-1/q\,'}a^{-1}(t)\|_{q\,'\!, \,(1,\infty)}\\
&\approx \|t^{-1/q}a(t)J(u(t),t;X_0,X_1)\|_{q,(1,\infty)}\\
& \le 2\, \|f\|_Y,
\end{align*}
i.e.,
\begin{equation}\label{f_infty_funkce2}
\|f_\infty\|_{X_0 \cap X_1}\lesssim  \|f\|_Y,
\end{equation}
which shows that \eqref{f_infty_funkce1} holds.

Now we are going to estimate $\|f\|_Z$. 
By \eqref{Afce111},
$$A(x)=a(x)\qquad \text{for all \ }x \in (0, 1).$$ 
Since $w(t)=0$ for all $t\ge 1$, we see that $J(w(t),t;X_0,X_0 \cap X_1)=0$ for all 
$t\ge 1$. Consequently,
\begin{equation}\label{redukcenormy}
 \hskip-0,1cm\  \|t^{-1/q}A(t)J(w(t),t;X_0,X_0 \cap X_1)\|_{q,(0,\infty)}=\|t^{-1/q}a(t)J(w(t),t;X_0,X_0 \cap X_1)\|_{q,(0, 1)}.
\end{equation}
Moreover, for all $t \in (0, 1)$,
$$
J(w(t),t;X_0,X_0 \cap X_1)=\max\left\{\|w(t)\|_{X_0},t\,\|w(t)\|_{X_0 \cap X_1}\right\},
$$
$$
\|w(t)\|_{X_0}\le \|u(t)\|_{X_0} +t\, \|f_\infty\|_{X_0}\le \|u(t)\|_{X_0} +t\, \|f_\infty\|_{X_0\cap X_1},
$$
$$
\|w(t)\|_{X_0\cap X_1} \le \|u(t)\|_{X_0\cap X_1}+t\, \|f_\infty\|_{X_0\cap X_1}.
$$
 Consequently, for all $t \in (0, 1)$,
\begin{align*}
&J(w(t),t;X_0,X_0 \cap X_1)\\
&\qquad \le \max\left\{\|u(t)\|_{X_0} +t\, \|f_\infty\|_{X_0\cap X_1}, 
t \left(\|u(t)\|_{X_0\cap X_1}+t\, \|f_\infty\|_{X_0\cap X_1}\right)\right\}\\
&\qquad\le \max\left\{\|u(t)\|_{X_0} +t\, \|f_\infty\|_{X_0\cap X_1}, 
t\,\|u(t)\|_{X_0\cap X_1}+t\, \|f_\infty\|_{X_0\cap X_1}\right\}\\
&\qquad=t\, \|f_\infty\|_{X_0\cap X_1}+\max\left\{\|u(t)\|_{X_0}, t\,\|u(t)\|_{X_0\cap X_1}\right\}\\
&\qquad=t\, \|f_\infty\|_{X_0\cap X_1}+J(u(t),t;X_0,X_0 \cap X_1).
\end{align*}
Using also \eqref{f_infty_funkce2} and the fact that (cf. Lemma \ref{1 o J funct})
$$
J(g,t;X_0,X_0 \cap X_1)=J(g,t;X_0, X_1) \quad \text{for all \ } g \in X_0 \cap X_1 \ \text{and \ } t\in(0, 1),
$$
we arrive at
\begin{equation}\label{odhad J_funk}
J(w(t),t;X_0,X_0 \cap X_1)\lesssim t\, \|f\|_Y+J(u(t),t;X_0, X_1)\quad \text{for all \ } t\in (0, 1).
\end{equation}
This, together with 
properties of slowly varying functions, implies that
\begin{align*}
{\rm RHS\eqref{redukcenormy}}&\lesssim \|t^{1-1/q}a(t)\|_{q\, (0, 1)} \,\|f\|_Y + 
\|t^{-1/q}a(t)J(u(t),t;X_0, X_1)\|_{q\, (0, 1)}   \\
&\approx \|f\|_Y + \|t^{-1/q}a(t)J(u(t),t;X_0, X_1)\|_{q\, (0, 1)}.  
\end{align*}
Thus, using also \eqref{nernorm}, we obtain that 
${\rm RHS\eqref{redukcenormy}}\lesssim \|f\|_Y.$
Together with \eqref{redukcenormy}, this yields
$$
\|t^{-1/q} A(t)J(w(t),t;X_0,X_0 \cap X_1)\|_{q,(0,\infty)}\lesssim \|f\|_Y.
$$
Finally, taking the infimum over all representations of the function $f$ in 
form \eqref{w1}, with $w:(0,\infty)\rightarrow X_0 \cap X_1$, we get 
\begin{equation}\label{Jembed}
Y=(X_0, X_1)_{0, q, a; J}\hookrightarrow Z=(X_0, X_0\cap X_1)_{0, q,A; J}.
\end{equation}

{\bf (ii)} Let $f \in Z$. Then there is a function $u:(0, \infty)\rightarrow X_0 \cap X_1$ such that 
\begin{equation}\label{f_funkceZ}
f=\int_0^\infty u(t) \frac{dt}{t}   \quad  {\rm ( {\text convergence\  in \ } X_0 )}
\end{equation}
and
\begin{equation}\label{nernormZ}
\|t^{-1/q}A(t)J(u(t),t;X_0,X_0\cap X_1)\|_{q,(0,\infty)}\le 2\,\|f\|_Z.
\end{equation}
Defining $f_{\infty}$ and $w$ by \eqref{f_infty_funkce} and \eqref{w}, we see that 
\eqref{w1} holds. We claim that \eqref{w v  pruniku} remains true. To verify this, it is again sufficient to show that 
\eqref{f_infty_funkce1} is satisfied. Using \eqref{f_infty_funkce}, we get
$$
\|f_\infty\|_{X_0 \cap X_1} \le \int_1^\infty \|u(t)\|_{X_0 \cap X_1}\frac{dt}{t}\,.
$$
Since 
$$
\|u(t)\|_{X_0 \cap X_1} \le \frac{J(u(t),t;X_0,X_0 \cap X_1)}{t}\quad \text{for all } t>0,
$$
we see that
$$
\|f_\infty\|_{X_0 \cap X_1} \lesssim \int_1^\infty \frac{J(u(t),t;X_0,X_0 \cap X_1)}{t}\frac{dt}{t}\,.
$$
Thus, applying H\" older's inequality, 
properties of slowly varying functions, and \eqref{nernormZ}, we obtain 
\begin{align}\label{new}
\|f_\infty\|_{X_0 \cap X_1} &\lesssim \|t^{-1/q}A(t)J(u(t),t;X_0,X_0 \cap X_1)\|_{q,(1,\infty)} \,\|t^{-1-1/q\,'}A^{-1}(t)\|_{q\,'\!,\, (1,\infty)}\\
&\approx \|t^{-1/q}A(t)J(u(t),t;X_0,X_0 \cap X_1)\|_{q,(1,\infty)}\notag\\
& \le 2\, \|f\|_Z,\notag
\end{align}
which shows that \eqref{f_infty_funkce1} holds.
Since 
$w(t)=0$ for all $t\ge 1$, we see that 
\begin{equation}\label{redukcenormy1}
  \|t^{-1/q}a(t)J(w(t),t;X_0,X_1)\|_{q,(0,\infty)}=\|t^{-1/q}a(t)J(w(t),t;X_0,X_1)\|_{q,(0, 1)}.
\end{equation}
Moreover, for all $t \in (0, 1)$,
$$
J(w(t),t;X_0,X_1)=\max\left\{\|w(t)\|_{X_0},t\,\|w(t)\|_{X_1}\right\},
$$
$$
\|w(t)\|_{X_0}\le \|u(t)\|_{X_0} +t\, \|f_\infty\|_{X_0}\le \|u(t)\|_{X_0} +t\, \|f_\infty\|_{X_0\cap X_1},
$$
$$
\|w(t)\|_{X_1} \le \|u(t)\|_{X_1}+t\, \|f_\infty\|_{X_1}
\le \|u(t)\|_{X_0\cap X_1} +t\, \|f_\infty\|_{X_0\cap X_1}.
$$
 Consequently, for all $t \in (0, 1)$,
\begin{align*}
J(w(t),t;X_0,X_1)
& \le \max\left\{\|u(t)\|_{X_0} +t\, \|f_\infty\|_{X_0\cap X_1}, 
t \left(\|u(t)\|_{X_0\cap X_1}+t\, \|f_\infty\|_{X_0\cap X_1}\right)\right\}\\
&\le t\, \|f_\infty\|_{X_0\cap X_1}+\max\left\{\|u(t)\|_{X_0}, t\,\|u(t)\|_{X_0\cap X_1}\right\}\\
&=t\, \|f_\infty\|_{X_0\cap X_1}+J(u(t),t;X_0,X_0 \cap X_1).
\end{align*}
Using also \eqref{new} 
we arrive at
\begin{equation}\label{2odhad J_funk}
J(w(t),t;X_0,X_1)\lesssim t\, \|f\|_Z+J(u(t),t;X_0,X_0 \cap X_1)\quad \text{for all \ } t\in (0, 1).
\end{equation}
Together with the fact that $a(x)=A(x)$ if $x\in (0,1)$, properties of slowly varying functions, and \eqref{nernormZ}, this implies 
that 
\begin{align*}
{\rm RHS\eqref{redukcenormy1}}&\lesssim  \|t^{1-1/q}a(t)\|_{q,(0,1)}\|f\|_Z+\|t^{-1/q}A(t)J(u(t),t;X_0,X_0 \cap  X_1)\|_{q\, (0, 1)}\\
&\lesssim \|f\|_Z.
\end{align*}
This estimate and \eqref{redukcenormy1} yield
$$
\|t^{-1/q}a(t)J(w(t),t;X_0,X_1)\|_{q,(0,\infty)}\lesssim \|f\|_Z.
$$
Finally, taking the infimum over all representations of the function $f$ in 
form \eqref{w1}, with $w:(0,\infty)\rightarrow X_0 \cap X_1$, we get 
$$
Z=(X_0, X_0\cap X_1)_{0, q,A; J}\hookrightarrow Y=(X_0, X_1)_{0, q, a; J},
$$
which, together with \eqref{Jembed}, shows that 
\begin{equation}\label{1equality}
(X_0, X_1)_{0, q, a; J}=(X_0, X_0\cap X_1)_{0, q,A; J}.
\end{equation}

{\bf (iii)} The function $A$ satisfies 
\begin{equation}\label{Aprop}
\int_0^x t^{-1}A^{-q'}(t)\, dt <\infty  \quad\text{for all \ }x>0, \qquad 
\int_0^\infty t^{-1}A^{-q'}(t)\, dt =\infty.
\end{equation}
Therefore, applying Theorem \ref{equivalence theorem} (with $(X_0, X_0\cap X_1), A$, and $B$ instead of $(X_0, X_1), a$, and $b$, 
respectively), we obtain 
$$ 
(X_0, X_0\cap X_1)_{0, q,A; J}=(X_0, X_0\cap X_1)_{0, q,B; K},
$$
and the proof of Theorem \ref{JS11} is complete.
\hskip 8cm$\square$

\medskip

{\bf II.} Proof of Theorem \ref{density4}. The function $A$ defined in Theorem \ref{JS11} satisfies \eqref{Aprop}. Thus, applying 
Theorem \ref{density theorem} (with $(X_0, X_1)$ replaced by $(X_0, X_0 \cap X_1)$ and with $a$ replaced by $A$) we obtain that 
the space $X_0 \cap X_1$ is dense in $(X_0, X_0\cap X_1)_{0,q,A;J}$. Together with \eqref{dual}, this yields the result. 
\hskip 9cm$\square$

\medskip
{\bf III.} Proof of Remark \ref{remark100J}.

{\bf (i)} The proof of \eqref{101J}. Put 
$$
Y:=(X_0,X_1)_{0,q,a;J}  \quad  \text{and \ }   \quad  Y_{(0,1)}:=(X_0,X_1)_{0,q,a;J;(0,1)}.
$$

Let $f \in Y$. Then there is a function $u:(0, \infty)\rightarrow X_0 \cap X_1$ such that 
\eqref{f_funkce} and \eqref{nernorm} hold. Defining the functions $f_{\infty}$ by 
\eqref{f_infty_funkce} and $w$ by 
\begin{equation}\label{wJ}
w(t):= u(t)+t\,f_\infty, \quad t \in (0, 1),
\end{equation}
we see that  
\begin{equation}\label{w11}
f=\int_0^1 w(t) \frac{dt}{t}.
\end{equation}
Moreover, we claim that
\begin{equation}\label{w v  prunikuJ}
w(t) \in  X_0 \cap X_1\quad \text{for all\ } t \in (0,1).
\end{equation}
 With respect to \eqref{wJ}, to verify \eqref{w v  prunikuJ}, it is sufficient to show that 
\eqref{f_infty_funkce1} is true. As in part~(i) of the proof of Theorem \ref{JS11}, we obtain 
\eqref{f_infty_funkce2}, which implies \eqref{f_infty_funkce1}.

Now we are going to estimate $\|f\|_{Y_{(0,1)}}$. Since, for all $t \in (0, 1)$,
$$
J(w(t),t;X_0,X_1)=\max\left\{\|w(t)\|_{X_0},t\,\|w(t)\|_{X_1}\right\},
$$
$$
\|w(t)\|_{X_0}\le \|u(t)\|_{X_0} +t\, \|f_\infty\|_{X_0}\le \|u(t)\|_{X_0} +t\, \|f_\infty\|_{X_0\cap X_1},
$$
$$
\|w(t)\|_{X_1} \le \|u(t)\|_{X_1}+t\,\|f_\infty\|_{X_1} \le \|u(t)\|_{X_1}+t\|f_\infty\|_{X_0\cap X_1},
$$
we obtain 
\begin{align*}
&J(w(t),t;X_0, X_1)\\
&\qquad \le \max\left\{\|u(t)\|_{X_0} +t\, \|f_\infty\|_{X_0\cap X_1}, 
t \left(\|u(t)\|_{X_1}+t\, \|f_\infty\|_{X_0\cap X_1}\right)\right\}\\
&\qquad\le \max\left\{\|u(t)\|_{X_0} +t\, \|f_\infty\|_{X_0\cap X_1}, 
t\,\|u(t)\|_{X_1}+t\, \|f_\infty\|_{X_0\cap X_1}\right\}\\
&\qquad=t\, \|f_\infty\|_{X_0\cap X_1}+\max\left\{\|u(t)\|_{X_0}, t\,\|u(t)\|_{X_1}\right\}\\
&\qquad=t\, \|f_\infty\|_{X_0\cap X_1}+J(u(t),t;X_0,X_1)\quad \text{for all \ }t \in (0, 1).
\end{align*}
Using \eqref{f_infty_funkce2}, we arrive at
$$
J(w(t),t;X_0,X_1)\lesssim t\, \|f\|_Y+J(u(t),t;X_0, X_1)\quad \text{for all \ } t\in (0, 1).
$$
Thus, applying also properties of slowly varying functions, we get
\begin{align*}
\|t^{-1/q}a(t)&J(w(t),t;X_0,X_1)\|_{q,(0, 1)}\\
&\lesssim \|t^{1-1/q}a(t)\|_{q\, (0, 1)} \,\|f\|_Y + 
\|t^{-1/q}a(t)J(u(t),t;X_0, X_1)\|_{q\, (0, 1)}   \\
&\approx \|f\|_Y + \|t^{-1/q}a(t)J(u(t),t;X_0, X_1)\|_{q\, (0, 1)},  
\end{align*}
which, together with \eqref{nernorm}, yields
$$
\|t^{-1/q} a(t)J(w(t),t;X_0,X_1)\|_{q,(0,1)}\lesssim \|f\|_Y.
$$
Finally, taking the infimum over all representations of the function $f$ in 
form \eqref{w11}, with $w:(0,1)\rightarrow X_0 \cap X_1$, we obtain 
$$
\|f\|_{Y_{(0,1)}}\lesssim \|f\|_Y\quad \text{for all \ } f\in Y, 
$$
i.e.,
\begin{equation}\label{Jembed0}
Y\hookrightarrow Y_{(0,1)}.
\end{equation}

Let $f \in Y_{(0,1)}$. Then there is a function $w:(0, 1)\rightarrow X_0 \cap X_1$ such that 
\eqref{w11} hold and 
\begin{equation}\label{0nernorm}
\|t^{-1/q}a(t)J(w(t),t;X_0,X_1)\|_{q,(0,1)}\le 2\,\|f\|_{Y_{(0,1)}}.
\end{equation}
Putting 
$$
u(t):=\begin{cases}w(t), &t \in (0, 1)\\
             0 \quad &t \in [1, \infty), 
		\end{cases}
$$
we see that  $u:(0,\infty)\rightarrow X_0 \cap X_1$,
\begin{equation}\label{Ef}
f=\int_0^\infty u(t)\, \frac{dt}{t}\,,
\end{equation}
and 
$$
J(u(t),t;X_0, X_1)=0\quad \text{for all \ } t\ge 1.
$$
Consequently,
\begin{align*}
\left\|t^{-1/q}\,a(t)\,J(u(t),t;X_0,X_1)\right\|_{q,(0,\infty)}
&=\left\|t^{-1/q}\,a(t)\,J(u(t),t;X_0,X_1)\right\|_{q,(0,1)}\\
&=\left\|t^{-1/q}\,a(t)\,J(w(t),t;X_0,X_1)\right\|_{q,(0,1)},
\end{align*}
which, together with \eqref{0nernorm}, yields
$$
\left\|t^{-1/q}\,a(t)\,J(u(t),t;X_0,X_1)\right\|_{q,(0,\infty)}\le  2\, \left\|f\right\|_{Y_{(0,1)}}.
$$
Thus, taking the infimum over all representations of the function $f$ in 
form \eqref{Ef}, with $u:(0,\infty)\rightarrow X_0 \cap X_1$, we get 
$$
 \|f\|_Y \lesssim\|f\|_{Y_{(0,1)}}   \quad \text{for all \ } f\in Y_{(0,1)}, 
$$
i.e.,
$$
 Y_{(0,1)}\hookrightarrow Y.
$$
This and \eqref{Jembed0} imply \eqref{101J}.

{\bf \,(ii)} The proof of \eqref{1002J} is quite analogous to that of \eqref{101J}.

{\bf (iii)} The proof of \eqref{1003J}. Put 
$$
Y:=(X_0,X_0\cap X_1)_{0,q,B;K}  \quad  \text{and \ }   \quad  Y_{(0,1)}:=(X_0,X_0\cap X_1)_{0,q,B;K;(0,1)}.
$$
The embedding $Y\hookrightarrow Y_{(0,1)}$ holds trivially. 

To prove the opposite embedding, assume that $f \in Y_{(0,1)}$. Let $Y_0:=X_0$ and $Y_1:=X_0\cap X_1$. 
Since $Y_1\hookrightarrow Y_0$, we obtain by Corollary \ref{col} (i) that
\begin{equation}\label{K1colY}
K(f,t;Y_0,Y_1)\approx\|f\|_{Y_0}\quad \text{for all} \ f \in Y_0 +Y_1=Y_0 \ \text{and}\ \,t \ge 1. 
\end{equation}
Note that the function $A$ satisfies \eqref{Aprop}. Thus, $\|t^{-1/q}B(t)\|_{q,(1,\infty)}<\infty$ 
(by Lemma \ref{102}~(ii) if $1<q<\infty$, and by \eqref{Bfce111} and \eqref{Aprop} if $q=\infty$). 
This and \eqref{K1colY} imply that, for all $f\in Y_0,$
\begin{align}\label{odhadK}
\|f\|_Y&=\|t^{-1/q}B(t)K(f,t;Y_0,Y_1)\|_{q,(0,\infty)}\\
&\approx \|t^{-1/q}B(t)K(f,t;Y_0,Y_1)\|_{q,(0,1)}+\|t^{-1/q}B(t)K(f,t;Y_0,Y_1)\|_{q,(1,\infty)}\notag\\
&\approx \|t^{-1/q}B(t)K(f,t;Y_0,Y_1)\|_{q,(0,1)}+\|t^{-1/q}B(t)\|_{q,(1,\infty)}\|f\|_{Y_0}\notag\\
&\approx \|t^{-1/q}B(t)K(f,t;Y_0,Y_1)\|_{q,(0,1)}+\|f\|_{Y_0}\notag\\
&= \|f\|_{Y_{(0,1)}}  +\|f\|_{Y_0}\notag.
\end{align}
Moreover, using \eqref{K1colY}, properties of the K-functional and slowly varying functions, we obtain, for all $f\in Y_0$,
\begin{align*}
\|f\|_{Y_0}&\approx K(f,1;Y_0,Y_1)
\approx\frac{K(f,1;Y_0,Y_1)}{1} \,\|t^{1-1/q}B(t)\|_{q,(0,1)}\\
&\le \|t^{-1/q}B(t)Kf,t;Y_0,Y_1)\|_{q,(0,1)}=\|f\|_{Y_{(0,1)}}.\\
\end{align*}
Together with \eqref{odhadK}, this shows that ${Y_{(0,1)}}\hookrightarrow Y.$ Therefore, \eqref{1003J} 
is true.
\hskip 2cm$\square$

\medskip
{\bf IV.} Proof of Corollary \ref{cor2}. 
Making use of the definition of the norm in the space $(X_0, X_0\cap X_1)_{0, q, B; K}$ and \eqref{1003J}, we get
\begin{align}\label{nk}
\|f\|_{(X_0,X_0\cap X_1)_{0,q,B;K}}&:=\|t^{-1/q}B(t)K(f,t;X_0,X_0\cap X_1)\|_{q,(0,\infty)}\\
&\,\approx \|t^{-1/q}B(t)K(f,t;X_0,X_0\cap X_1)\|_{q,(0,1)}\notag
\end{align}
for all $f \in X_0+(X_0\cap X_1)=X_0.$ 
Moreover, by Lemma \ref{o K funct},
\begin{equation}\label{kner}
K(f,t;X_0,X_0\cap X_1)\approx t\,\|f\|_{X_0}+K(f,t;X_0,X_1)\quad \text{for all}\ \,f\in X_0 \ \text{and}\ t\in (0,1). 
\end{equation}
Thus, applying also properties of slowly varying functions, we get, for all $f \in  X_0$, 
\begin{align}\label{kner1}
\|f\|_{(X_0,X_0\cap X_1)_{0,q,B;K}}&\approx \|t^{1-1/q}B(t)\|_{q,(0,1)} \,\|f\|_{X_0}+\|t^{-1/q}B(t)K(f,t;X_0,X_1)\|_{q,(0,1)}\\
&\approx \|f\|_{X_0}+\|t^{-1/q}B(t)K(f,t;X_0,X_1)\|_{q,(0,1)}.\notag
\end{align}
The embedding $X_0\cap X_1 \hookrightarrow X_1$ implies that 
$$
K(f,t;X_0,X_1)\lesssim  K(f,t;X_0,X_0\cap X_1)\quad \text{for all}\ \,f\in X_0 \ \text{and}\ t>0. 
$$
Consequently, for all $f \in  X_0$,
\begin{align}\label{kner2}
\|t^{-1/q}B(t)K(f,t;X_0,X_1)\|_{q,(1,\infty)}&\lesssim \|t^{-1/q}B(t)K(f,t;X_0,X_0\cap X_1)\|_{q,(1,\infty)}\\
&\le \|t^{-1/q}B(t)K(f,t;X_0,X_0\cap X_1)\|_{q,(0,\infty)}\notag\\
&=\|f\|_{(X_0,X_0\cap X_1)_{0,q,B;K}}.\notag
\end{align}
Making use of \eqref{kner1}, \eqref{kner2}, and the definition of the norm in 
$(X_0,X_1)_{0,q,B;K}$, we arrive at
\begin{align*}
\|f\|_{X_0}+&\|f\|_{(X_0,X_1)_{0,q,B;K}}\\
&\approx \|f\|_{X_0} + \|t^{-1/q}B(t)K(f,t;X_0,X_1)\|_{q,(0,1)}+ \|t^{-1/q}B(t)K(f,t;X_0,X_1)\|_{q,(1,\infty)}\\
&\lesssim \|f\|_{(X_0,X_0\cap X_1)_{0,q,B;K}}
\lesssim \|f\|_{X_0}+\|f\|_{(X_0,X_1)_{0,q,B;K}}\quad \text{for all}\ \,f\in X_0,
\end{align*}
which implies that
\begin{equation}\label{kner3}
(X_0,X_0\cap X_1)_{0,q,B;K} = X_0 \cap (X_0,X_1)_{0,q,B;K}.
\end{equation}
Furthermore, the function $A$ satisfies \eqref{Aprop}. 
Consequently, by Theorem \ref{equivalence theorem} (with $a$ and $b$ replaced by $A$ and $B$, respectively),
\begin{equation}\label{rovnost1}
(X_0,X_1)_{0,q,B;K} = (X_0,X_1)_{0,q,A;J}.
\end{equation}
Finally, using \eqref{dual}, \eqref{kner3}, and \eqref{rovnost1}, we arrive at
$$
(X_0,X_1)_{0,q,a;J}=(X_0,X_0\cap X_1)_{0,q,B;K}= X_0 \cap (X_0,X_1)_{0,q,B;K}= X_0 \cap (X_0,X_1)_{0,q,A;J},
$$
and \eqref{corres2} is verified.
\hskip 11,2cm$\square$

\medskip

\section {Proofs of Theorem \ref{*equivalence theorem1*}, Remark \ref{remark100*}, and Corollary \ref{cor2/3}}\label{Pr1Main*}

{\bf I.} Proof of Theorem \ref{*equivalence theorem1*}. As in the proof of Theorem \ref{KS}, 
we obtain \eqref{31211} with $q=\infty$, i.e.,
\begin{equation}\label{31211*}
(X_0, X_1)_{0,\infty,b;K}=(X_0, X_0+X_1)_{0,\infty,B;K}.
\end{equation}
Moreover, the function $B$ satisfies $B\in \, SV(0,\infty)\cap AC(0,\infty),$ 
\begin{equation}\label{*B*}
B \text{\ is strictly decreasing},\quad B(0)=c\,\beta(0)=\infty, \quad B(\infty)=b(\infty)=0.\text{\footnotemark$^)$}
\end{equation}
\footnotetext{$^{)}$ Condition \eqref{*B*} implies that the analogue of \eqref{B*} with $q=\infty$ holds, i.e.,
$\|B(t)\|_{\infty,(x,\infty)}<\infty$ for all $x>0$, $\|B(t)\|_{\infty,(0,\infty)}=\infty$.}
Thus, applying Theorem \ref{equivalence theorem1*} (with $(X_0, X_1), b,$ and $a$ replaced by $(X_0, X_0+X_1), B,$  and $A$, 
respectively), we arrive at 
\begin{equation}\label{B*1*}
(X_0, X_0+X_1)_{0,\infty,B;K}=(X_0, X_0+X_1)_{0,\infty,A;J}.
\end{equation}
Together with \eqref{31211*}, this gives \eqref{D2E1*}.
\hskip 8,4cm$\square$

\medskip
{\bf II.} Proof of Remark \ref{remark100*}. The proofs of \eqref{101*} and \eqref{1002*} are quite analogous to those of 
\eqref{101} and \eqref{1002}.
 
Moreover, if $\|A(t)\|_{\infty, (1,e)}<\infty$, then also the proof of \eqref{1003*} can be done analogously to that of 
\eqref{1003}. (Note that condition \eqref{3003} with $q=\infty$ can be verified on making use of \eqref{difeq}.)
\hskip 13,45cm$\square$

\medskip
{\bf III.} Proof of Corollary \ref{cor2/3}. The proof is analogous to that of Corollary \ref{cor2/2} (the only difference is that 
we apply Theorem \ref{equivalence theorem1*} instead of Theorem \ref{equivalence theorem1}.)
\hskip 4,4cm$\square$

\medskip
\section{Proofs of Theorems \ref{equivalence theorem1**a}, \ref{density theorem**}, Remark \ref{remark100J1}, and Corollary \ref{cor21}}\label{Pr1MainJ1}

{\bf I.} Proof of Theorem \ref{equivalence theorem1**a}. As in the proof of Theorem \ref{JS11}, 
we obtain \eqref{1equality} with $q=1$, i.e.,
\begin{equation}\label{1equality1}
(X_0, X_1)_{0, 1, a; J}=(X_0, X_0\cap X_1)_{0, 1,A; J}.
\end{equation}
Moreover, the function $A$ satisfies $A\in \, SV(0,\infty)\cap AC(0,\infty),$ 
\begin{equation}\label{*A*}
A \text{\ is strictly decreasing},\quad A(0)=a(0)=\infty,\quad A(\infty)=c\,\alpha(\infty)=0.\text{\footnotemark$^)$}
\end{equation}
\footnotetext{$^{)}$ Condition \eqref{*A*} implies that the analogue of \eqref{Aprop} with $q=1$ holds, i.e.,
$\|A^{-1}(t)\|_{\infty,(0,x)}<\infty$ for all $x>0$, $\|A^{-1}\|_{\infty,(0,\infty)}=\infty$.}
Thus, applying Theorem \ref{equivalence theorem1**} (with $(X_0, X_1), a,$ and $b$ replaced by $(X_0, X_0\cap X_1), A,$  and $B$, 
respectively), we arrive at 
\begin{equation}\label{B*1*1}
(X_0, X_0\cap X_1)_{0,1,A;J}=(X_0, X_0\cap X_1)_{0,1,B;K}.
\end{equation}
Together with \eqref{1equality1}, this gives \eqref{D2E1*S}.
\hskip 8,1cm$\square$

\medskip
{\bf II.} Proof of Theorem \ref{density theorem**}. The function $A$ defined in Theorem \ref{equivalence theorem1**a}
satisfies \eqref{*A*}. Thus, applying Theorem \ref{density theorem*} (with $(X_0, X_1)$ and $a$ replaced by $(X_0, X_0\cap X_1)$  and $A$, 
respectively), we get the result.
\hskip 11,7cm$\square$

\medskip
{\bf III.} Proof of Remark \ref{remark100J1}. The Proofs of \eqref{101J1} and \eqref{1002J1} are analogous to that of \eqref{101J}.
 Since 
$$
\|A^{-1}(t)\|_{\infty,(0,x)}<\infty \quad \text{for all }x>0, \qquad 
\|A^{-1}(t)\|_{\infty,(0,\infty)}=\infty
$$
(which is an analogue of \eqref{Aprop} with $q=1$),  $\|t^{-1}B(t)\|_{1, (1,\infty)}<\infty$, and
$\|B(t)\|_{1,(0,1)}<\infty$ (which means that $\|t^{1-1/q}B(t)\|_{q,(0,1)}<\infty$ holds with $q=1$),
one can prove \eqref{1003J1} similarly to \eqref{1003J}. 
\hskip 13,1cm$\square$

\medskip
{\bf IV.} Proof of Corollary \ref{cor21}. We proceed analogously as in the proof of Corollary \ref{cor2}. Since 
$\|B(t)\|_{1,(0,1)}<\infty$, we see that \eqref{kner1} with $q=1$ is true. Also \eqref{kner2} and \eqref{kner3} 
with $q=1$ hold. Furthermore, the function $A$ satisfies \eqref{*A*}. Consequently, applying 
Theorem~\ref{equivalence theorem1**}  (with $A$ and $B$ instead of $a$ and $b$, respectively),
 we get \eqref{rovnost1} with $q=1$. Thus, making use of \eqref{D2E1*S}, further \eqref{kner3} and \eqref{rovnost1} (both with $q=1$), we obtain 
$$
(X_0,X_1)_{0,1,a;J}=(X_0,X_0\cap X_1)_{0,1,B;K}= X_0 \cap (X_0,X_1)_{0,1,B;K}= X_0 \cap (X_0,X_1)_{0,1,A;J},
$$
and \eqref{corres21} is verified.
\hskip 11,0cm$\square$

\bibliographystyle{alpha}

\def\cprime{$'$}

\end{document}